\documentclass[]{amsart}
\usepackage[]{amsmath, amsthm, amsfonts, amssymb }
\usepackage{latexsym,amsxtra,amscd,ifthen}
\usepackage{verbatim}
\usepackage[all]{xy}

\newcommand {\red} {{\textrm red}}

\newcommand {\N} {{\mathbb N}}
\newcommand {\C} {{\mathbb C}}
\newcommand {\R} {{\mathbb R}}
\newcommand {\Z} {{\mathbb Z}}
\newcommand {\Q} {{\mathbb Q}}
\newcommand {\PP} {{\mathbb P}}

\newcommand {\F} {{\mathcal F}}

\newcommand {\BS} {{\mathbb S}}
\newcommand {\E} {{\mathcal E}}
\newcommand {\dt} {{\bullet}}
\newcommand {\calC} {{\mathcal C}}
\newcommand {\G} {{\mathcal G}}

\newcommand {\calN} {{\mathcal N}}
\newcommand {\calI} {{\mathcal I}}
\newcommand {\cp} {{$F$-ample }}

\newtheorem{thm}{Theorem}
\newtheorem{cor}[subsection]{Corollary}
\newtheorem{lemma}[subsection]{Lemma}
\newtheorem{prop}[subsection]{Proposition}
\newtheorem{defn}{Definition}
\newtheorem{remark}[subsection]{Remark}
\newtheorem{ex}[subsection]{Example}

\begin{document}

\title[Frobenius amplitude]{
Frobenius amplitude and strong vanishing theorems for
vector bundles
}
\author{
        Donu Arapura    
}
\address{Arapura: Department of Mathematics\\
Purdue University\\
West Lafayette, IN 47907}
\email{dvb@math.purdue.edu}
\thanks{Arapura partially supported by the NSF}
\address{Keeler: Department of Mathematics \\ MIT \\ Cambridge, MA
  02139-4307   }
  \thanks{ 
     Keeler partially supported by an NSF Postdoctoral Research Fellowship.}
\email{dskeeler@mit.edu}
\urladdr{http://www.mit.edu/\~{}dskeeler}


\date{} 

\maketitle

\begin{center} With Appendices by Dennis S. Keeler \end{center}

\begin{abstract}
The primary goal of this paper is to systematically exploit the method
of Deligne-Illusie to obtain Kodaira type vanishing theorems for
vector bundles and more generally coherent sheaves on algebraic
varieties. The key idea is to introduce a number which provides a
cohomological measure of the positivity of a coherent sheaf called the
Frobenius or F-amplitude.
The F-amplitude enters into the statement of the basic vanishing
theorem, and this leads to the problem of calculating, or at least
estimating, this number. Most of the work in this paper is devoted to
doing this various situations.  
\end{abstract}

\tableofcontents

In \cite{deligne-ill}, Deligne, Illusie and Raynaud gave a beautiful
 proof of the Kodaira-Akizuki-Nakano vanishing theorem
for ample line bundles using characteristic $p$ methods.
The  goal of this paper is to apply these methods to  
obtain vanishing theorems for vector bundles and, more generally, 
 sheaves in a systematic fashion. 
In order to facilitate this, we  introduce a cohomological measure of
the positivity of  a coherent sheaf on an algebraic variety that we call the
Frobenius or $F$-amplitude. 
We also introduce some variations on this idea, such as  the 
$F$-amplitude  relative to a normal crossing divisor.
The smaller the amplitude, the more
positive it is; when it is zero, we say that the sheaf is $F$-ample.
($F$-ample vector bundles have been called ``cohomologically $p$-ample''
in  \cite{gieseker}, \cite{mi} and possibly elsewhere, but we prefer
 the shorter term.) $F$-ampleness for  bundles of
rank greater than one turns out to be an unreasonably restrictive notion, and it
appears more useful to consider the class of bundles with 
small  $F$-amplitude relative to the rank.

As the terminology suggests, the definition of $F$-amplitude
makes use of the Frobenius map  in an essential way.
However, it can be extended into characteristic zero by the usual 
reduction modulo $p$ tricks. While this leads to a definition,
it is one that is not particularly convenient to use in practice.
For curves and projective spaces, we can give a 
reformulation of  $F$-amplitude in  characteristic free terms.
In general, it seems that the best
we can hope for are some reasonable  bounds on $F$-amplitude, and
much of this paper is devoted to finding such bounds. 
The key result in this direction is 
 theorem~\ref{thm:keyestimate}, which shows that in characteristic 
 zero the $F$-amplitude of an ample vector bundle is bounded above by its
 rank. The proof relies on some work of Carter and Lusztig in
 modular representation theory.

 The penultimate section contains the main theorem. It gives
the vanishing of the cohomology groups of a sheaf on a smooth projective
variety tensored with the differentials with logarithmic singularities
along a divisor in a range determined by
the $F$-amplitude relative to the divisor. 
A special case of this for $F$-ample bundles had been considered by
Migliorini \cite{mi}. The vanishing theorem is 
nominally a characteristic $p$ result; the interesting 
consequences are in characteristic zero. From this we are able to
 recover some old results such as Le Potier's vanishing theorem, 
and to discover some new ones as well.
One corollary that we want to call attention to is the following 
Kawamata-Viehweg type theorem (cor.~\ref{cor:kvvan1}):
Let $\E$ be a vector bundle on smooth projective variety $X$. Suppose
there is an effective fractional $\Q$-divisor $\Delta$ with normal crossing 
support $D$ such that $\E(-\Delta)$ is ample, which means that some 
symmetric power $S^m(\E)(-m\Delta)$ is ample in the usual sense.
Then $H^i(\Omega_X^j(log D)(-D)\otimes \E) = 0$ for $i+j\ge dim\,X+ rank(\E)$;
in particular, $H^i(\omega_X \otimes \E) = 0$ for $i\ge rank(\E)$.
This result  is put to use in the final section
to obtain a refinement of the Lefschetz
hyperplane theorem and to obtain a Le Potier
theorem for noncompact varieties.

 The notion of $F$-semipositivity is obtained by relaxing the
 condition  for $F$-ampleness. We show that $F$-semipositive vector bundles
are nef. In characteristic $0$, more is true, namely  $F$-semipositive  bundles
are ``arithmetically nef'' which  means roughly that it specializes 
to a nef bundle in positive characteristic. The converse fails in general.
However,  for line bundles the equivalence of these notions has been established
by  Dennis Keeler, and included as an appendix.
This can be used to slightly extend the aforementioned vanishing theorem.

\section{Frobenius amplitude}\label{sec:frobenius-amplitude}

In this section, we define the notion
of Frobenius (or simply $F$-) amplitude.
This definition is most natural in positive characteristic,
and we start with this case.
Let $k$ be a field of characteristic $p > 0$, and
let $X$ be a  variety defined over $k$.
 $F$, or sometimes $F_X$, will denote the absolute Frobenius of 
$X$ (i.e. the morphism of schemes which is the
identity on the set $X$ and the $p$th power map on $O_X$).
The absolute Frobenius can be factored as:
$$
\xymatrix{
 X\ar[r]_{F'}\ar[rd]\ar@/^/[rr]^{F_X} & X'\ar[d]\ar[r]_{F_k} & X\ar[d] \\
  & spec\, k\ar[r]_{F_k} & spec\, k
}
$$
where the righthand square is cartesian. $F'$ is the relative
Frobenius. When $k$ is perfect, $F_k:spec\,k\to spec\,k$ and its base
change $X'\to X$ are isomorphisms of $\Z/p\Z$-schemes. In view of this,
the relative Frobenius can be replaced by the absolute
Frobenius and $X'$ by $X$ in the statements of  
\cite[2.1, 4.2]{deligne-ill}.

Given a coherent sheaf $\E$, denote $F^{n*}\E$ by
$\E^{(p^n)}$. For a vector bundle $\E$ given
by a $1$-cocycle $g_{ij}$, $\E^{(p^n)}$
is given by  $g_{ij}^{p^n}$. If $\mathcal{I}$ is an ideal sheaf on $\PP^n$
generated by polynomials $f_i$, then $\mathcal{I}^{(p^n)}$ is the ideal
sheaf generated by $f_i^{p^n}$.
Define  the  {\em  $F$-amplitude} 
$\phi(\E)$ of a coherent sheaf $\E$ 
to be  the smallest integer $l$ such that
 for any locally free sheaf $\F$, there
exists an $N$ such that $H^i(X, \E^{(p^m)}\otimes \F) = 0$ 
for all $i > l$ and $m > N$.
A few words of caution should be added here. We are purposely using the
naive definition, but this  has reasonable properties only
when $X$ is smooth (which implies that $F$ is flat) or $\E$ is locally
free. In more general situations, $F^{n*}\E$ should be replaced by the
derived pullback $LF^{n*}\E$, at which point $\E$ may as well be replaced by
an object in $D_{coh}(X)$ (one day, perhaps).
We have that $\phi(\E)$ is less than or equal to the coherent
cohomological dimension of $X$ which is less than
or equal to the $dim X$.

Now suppose that $k$ is a field of characteristic $0$.
By a diagram over a scheme $S$, we will mean a collection of
$S$-schemes $X_i$, $S$-scheme morphisms $f_{ij}:X_{i}\to X_{j}$,
$O_{X_{i}}$-modules $\E_{i,l}$ and morphisms between the 
pullbacks and pushforwards of these modules.
Given a morphism $S'\to S$,  and a diagram $D$ over $S$, we
can define its fiber product $D\times_{S}S'$ in the 
obvious way.
Given a diagram $D$ over $Spec\,k$,
an {\em arithmetic thickening} (or simply just thickening)  of it is
a choice  of a finitely generated
$\Z$-subalgebra $A\subset k$, and a diagram $\tilde D$ over $Spec\, A$,
so that $D$ is isomorphic to the fiber product over $Spec\, k$.
Given two thickenings $\tilde D_{i}\to Spec\, A_{i}$,
we will say the second refines the first if there is
a homomorphism $A_1\to A_2$, and an isomorphism
between $D_2$ and $D_1\times_{Spec A_1}Spec A_2$.

By standard arguments (e. g.  \cite[sect. 6]{illusie2}):

\begin{lemma}
Any finite diagram of $k$-schemes of finite type and
 coherent sheaves has an arithmetic thickening.
 Any two thickenings have a common refinement.
\end{lemma} 

Suppose that $X$ is a  quasiprojective
$k$-variety with a coherent sheaf $\E$.
Given a thickening $(\tilde X, \tilde \E)$
over $A$,
we will write $p(q) = char(A/q)$, $X_q$ for the fiber and 
$\E_q = \tilde\E|_{X_{q}}$ for each closed point $q\in Spec A$. 
We will say that a property holds for {\em almost all $q$} if
it holds for all $q$ in a nonempty open subset of $Spec\, A$.
For each closed point $q\in Spec\, A$, the fiber $X_q$
is defined over the finite field $A/q$, so
that the $F$-amplitude of the restriction $\E_q$
can be defined as before. 
We say that $i\ge\phi(\E)$ if and only if $i\ge\phi(\E_q)$ holds for
almost all $q$. Equivalently, 
the  $F$-amplitude $\phi(\E)$ is obtained by minimizing
$\max_q\,\phi(\E_q)$  over all thickenings. 
Note that there is no (obvious) semicontinuity property for
$\phi(\E_q)$. So it is not clear if this is the optimal
definition, but it is sufficient for the present purposes.
Any alternative definition should satisfy the following:
for any arithmetic thickening of a finite 
collection of coherent sheaves $\E_1,\ldots \E_N$, there is a 
 sequence of closed points $q_j$ with $char\, (A/q_j)\to \infty$
such that  $\phi(\E_i)\ge \phi((\E_i)_{q_j})$.

Let $X$ be a smooth projective variety over 
$k$. We have an ordering on divisors defined in the usual
way: $D\le D'$ if and only if
the coefficients of $D$ are less than or equal to the
coefficients of $D'$. 
Fix  a reduced divisor $D\subset X$ with normal
crossings. 
Assume that $char\, k =p >0$, then we define
the $F$-amplitude of a coherent sheaf $\E$ relative to $D$
 as follows
$$\phi(\E,D) = \min\{\phi(\E^{(p^n)}(-D'))\, |\,n\in \N,\, 0\le 
D'\le (p^n-1)D\}. $$
If $D'\le (p^n-1)D$ is a divisor for which this minimum is achieved, we will
refer to the $\Q$-divisor $\frac{1}{p^n}D'$ as a critical divisor for
$\E$ relative to $D$.
It will be  convenient to introduce the relation on divisors,
$A<_{strict} B$ if the multiplicity of $A$ along any irreducible
component $C$ of the union of their supports
 is less than the multiplicity of $B$ along $C$.
Then the upper inequality above is just that $D'<_{strict} p^nD$.
When $char\, k = 0$, we proceed as above, $\phi(\E,D)$
is the minimum of $\max_q\, \phi(\E_q,D_q)$ over all thickenings
of $(X,D, \E)$. We define the generic $F$-amplitude
$\phi_{gen}(\E)$ of a locally free sheaf $\E$ to be the infimum of
$\phi(f^*\E, D)$ where $f:Y\to X$ varies over all birational
maps $f$ with exceptional divisor $D$ such that
$Y$ is smooth and $D$ has normal crossings.

In any characteristic, we will define $\E$ to be \cp if
and only if $\phi(\E) = 0$. We will see below that
a line bundle is ample if and only if it is  $F$-ample.
However for bundles of higher rank, $F$-ampleness
is a stronger condition. In positive characteristic,
$F$-ample vector bundles are the same as cohomologically
$p$-ample vector bundles as defined in \cite{gieseker}.

Most of the work below will be in positive characteristic.
The proofs in characteristic zero are  handled
by standard semicontinuity arguments on a thickening.

Throughout the rest of this paper, unless stated otherwise,
$X$ will denote a  projective variety over a 
field $k$, and the  symbols $\E, \F,\ldots$ will denote coherent sheaves
on $X$.

\section{Elementary bounds on $F$-amplitude}

\begin{lemma}\label{lemma:one}
    If a sheaf $\F$ on a topological space is 
 quasi-isomorphic to a bounded complex $\F^\dt$
then  $H^i(\F) = 0$  provided that $H^a(\F^b) = 0$ for all
$a+b=i$.
\end{lemma}

\begin{proof} This follows
 from the spectral sequence
 $$E_1^{ab}=H^b(\F^a)\Rightarrow \mathbb{H}^{a+b}(\F^\dt)\cong H^{a+b}(\F).$$
 \end{proof}
    
\begin{lemma}\label{lemma:locfree}
Suppose  $char\, k =p > 0$ and that $\E$ is a locally free sheaf
on $X$. Then for any coherent sheaf $\F$,
    $$H^i(X, \E^{(p^m)}\otimes \F) = 0$$ for $i>\phi(\E)$ and $m>>0$.
\end{lemma}   

\begin{proof}
This will be proved by descending induction starting from
$i=dim\,X+1$.
Choose an ample line bundle $O(1)$.
 We can find an exact sequence
$$0\to \F''\to \F'\to \F\to 0$$
where $\F'$ is a sum of twists of $O(1)$  
(by Serre's theorems, we can take $\F'= H^0(\F(n))\otimes O(-n)$ for $n>>0$).
Tensoring this with $\E^{(p^m)}$ and applying the long
exact sequence for cohomology shows that $H^i(X, \E^{(p^m)}\otimes \F) = 0$
for $i>\phi(\E)$ and $m>>0$.
\end{proof}

The proof gives something slightly stronger:

\begin{cor}\label{cor:FamplViaLineb}
Fix an ample line bundle $O_X(1)$. Then 
$\phi(\E)\le A$ if and only if for any $b$ there exists
$n_0$ such that 
$$H^i(\E^{(p^n)}(b)) = 0$$
for all $i> A$, $n\ge n_0$.
\end{cor}

\begin{lemma}\label{lemma:ampleline}
 A line bundle is $F$-ample if and only if it is ample.
\end{lemma}

\begin{proof}
First assume that we are over a field of characteristic $p>0$.    
Then we have $L^{(p^{n})} = L^{p^n}$. Therefore if $L$ is ample,
it is $F$-ample by Serre's vanishing theorem. Suppose
$L$ is $F$-ample. Choose $x_0\in X$, then by lemma \ref{lemma:locfree}
$H^1(X,m_{x_0}\otimes L^{n_0}) = 0$ for $x_0$ and some
$n_0$ a power of $p$. Therefore $L^{n_0}$ has a global
section $s_0$ which is nonzero at $x_0$.
Let $U_0$ be the complement of the zero set of $s_0$.
If $U_0\not= X$, we can choose $x_1$ in the complement and
arrange that $H^1(X,m_{x_1}\otimes L^{n_1}) = 0$ for 
some power $n_1$ of $p$. Therefore $L^{n_1}$ has a section
$s_1$ not vanishing at $x_1$.
If $U_0\cup U_1\not= X$, then we can choose $x_2$ in the complement
and proceed as above. Eventually this process has to stop, because
$X$ is noetherian. Therefore $L^{n_0n_1\ldots}$ is generated
by the sections $s_0^{n_1n_2\ldots}, s_1^{n_0n_2\ldots}, \ldots$.
Repeating the same line of reasoning with the sheaves
$ m_xm_y\otimes L^{n}$ shows that some power of $L$
is very ample.

Choose a thickening of $(\tilde X, \tilde L)$ over $A$.
If $L$ is ample, then we can assume $\tilde L$ is ample
by shrinking $Spec\, A$ if necessary. Consequently $\tilde L_{q}$
is $F$-ample for each closed point $q\in Spec\, A$ by the previous
paragraph. Therefore $L$ is $F$-ample.
Now suppose that $L$ is $F$-ample. As above, it suffices to
show that for any ideal sheaf $I$,
$H^1(I\otimes L^{n}) = 0$ for some $n >0$. But this is
easily seen by choosing a thickening of $(X,L,I)$
applying the previous case on a closed fiber, 
using semicontinuity to deduce this for the generic fiber,
then flat base change to deduce the vanishing for $X$.
\end{proof}

\begin{thm}\label{thm:first}
Let $\E, \E_0\ldots $ be
coherent sheaves on a projective variety $X$. Assume either that
$X$ is smooth or that these sheaves are locally free.
Then the following statements hold.
\begin{enumerate}

    \item\label{thm:first1} Given
     an exact sequence
     $0\to\E_1\to \E_2\to \E_3\to 0$,
    $\phi(\E_2) \le  max(\phi (\E_1), \phi (\E_3))$.

    \item\label{thm:first2} Let
    $$0\to \E_n\to \E_{n-1}\to \ldots \E_0\to \E\to 0$$
    be an exact sequence   such  that  for each $i$, 
   $\phi (\E_i)\le i +l$,
    then $\phi (\E) \le l$.

  \item\label{thm:first3}
 Let $0\to \E\to \E^0\to \E^1\to \ldots \E^n\to 0$
 be an exact sequence 
 such  that  for each $i$, $\phi (\E^i) \le l-i$
    then $\phi (\E) \le l$.

    \item\label{thm:first4}
    Let $f:Y\to X$ be a proper morphism of  projective varieties
     such that $d$ is the maximum dimension of the closed fibers. If
    $\E$ is locally free then
    $\phi(f^*\E) \le \phi(\E) + d$.
    In particular, if $f$ is a closed immersion, then $\phi(f^*\E)\le \phi(\E)$.
    
    \item\label{thm:first5} If $f:Y\to X$ is an \'etale
morphism of smooth projective varieties, then $\phi(f_*\E) = \phi(\E)$.

\end{enumerate}
\end{thm}

\begin{proof}
 The first statement is obvious, and second and third follow
 from lemma \ref{lemma:one}

 For the remaining statements, we will assume that $char k = p > 0$,
the characteristic $0$ case is a straightforward semicontinuity
argument.
There is a commutative diagram
\begin{equation}\label{eq:frobenius-diagram}
\xymatrix{
 Y\ar[d]\ar[r]^{F^m_Y} & Y\ar[d] \\
 X\ar[r]^{F^m_X} & X
}
\end{equation}
 Suppose that $f$ is proper with fibers  of dimension $\le d$.
  If $\E$ is a locally free $O_X$-module, and $\F$ a coherent
 $O_Y$-module, then 
 $ H^i(\E^{(p^m)}\otimes R^jf_*\F) = 0$ for $i>\phi(\E)$
 and $m >>0$.
 Therefore the Leray spectral sequence implies
 $$ H^{i}((f^*\E)^{(p^m)}\otimes \F)  =0$$
 for $i>\phi(\E)+d$ and $m >> 0$.
 
 If $f:Y\to X$ is \'etale, then the above diagram is cartesian.
Furthermore, $F_X$ and $F_Y$ are both flat when $X$ and $Y$
are smooth. Cohomology commutes with 
flat base change, therefore for any coherent $O_Y$-module
$\E$ and locally free $O_X$-module $\F$,
$$H^i(X, (f_*\E)^{(p^m)}\otimes \F) \cong 
H^i(Y, \E^{(p^m)}\otimes f^*\F),$$
and this implies the equality of amplitudes. 
\end{proof}

These result easily extend to the case
of $F$-amplititude relative to a divisor. Here we just
treat one case that will be needed later.

\begin{lemma}\label{lemma:genfinitepullback}
Let $Y\to X$ be a morphism of varieties with $Y$ smooth.
Suppose that $D= \sum D_i$ is a divisor with normal
crossings on $Y$, such that there exist $a_i\ge 0$
for which $L=O_Y(-\sum a_iD_i)$ is
relatively ample. Then for any locally free sheaf $\E$
on $X$, $\phi(f^*\E, D)\le \phi(\E)$.
\end{lemma}

\begin{proof}
The proof is very similar  case (4), of the previous theorem.
Assume $char\, k = p$, and choose $n_0$ such that
$p^{n_0}> a_i$. Set $\E' = f^*\E^{(p^{n_0})}\otimes L$
and let $\F$ be another coherent sheaf on $Y$.
Since $L$ is relatively ample, the higher direct images
of $\F\otimes L^N$ vanish for $N>>0$.
Therefore, the spectral sequence
$$H^a(\E^{(p^{n+n_0})}\otimes R^bf_*(\F\otimes L^{p^n}))
\Rightarrow H^{i}((\E')^{(p^n)}\otimes \F)$$
yields the vanishing of the abutment for $i>\phi(\E)$ and
$n>>0$.
\end{proof}

\begin{cor}\label{cor:biratphi}
Suppose that $k$ has characteristic $0$. 
If $f:Y\to X$ is a resolution of singularities such that the
exceptional divisor $D$ has normal crossings, then 
$\phi(f^*\E, D)\le \phi(\E)$.
If $g:Z\to X$ is any resolution of singularities,
then $\phi_{gen}(g^*\E)\le \phi(\E)$.
\end{cor}

\begin{proof} 
Since $f$ can be realized as the blow up of
$X$ along an ideal, it follows that we can find
a relatively ample divisor of the form $-\sum\, a_iD_i$.
with $a_i\ge 0$ where the $D_i$ are the irreducible
components of $D_{red}$. The first assertion clearly implies that
second,  since  $Z$ can be blown up further.
\end{proof}

\begin{cor}\label{cor:genphi}
Suppose that $k$ has characterisitic $0$. If $g:Z\to X$ is a surjective
morphism  with $Z$ smooth, then $\phi_{gen}(g^*\E)\le \phi(\E)+d$,
where $d$ is the dimension of the generic fiber.
\end{cor}

\begin{proof}
Construct the following commutative diagram:
$$
\xymatrix{
  & Y\ar[r]^{\alpha}\ar[d]^{\beta}\ar[ld]^{\epsilon} & Z\ar[d]^{g}\ar@{-->}[ld]^{\gamma} \\
 Y'\ar[r]_{\kappa} & \PP^d\times X\ar[r]_{\delta} & X
}
$$
where $\gamma$ is a generically
finite rational map (which exists by Noether's normalization lemma), 
 $\delta$ is the projection,
$\alpha$ is a resolution  of the indeterminacy locus of
$\gamma$, and $Y\to Y'\to \PP^d\times X$
the Stein factorization of $\beta$. By theorem \ref{thm:first} (4) and
the previous corollary (applied to $\delta\circ \kappa$ and $\epsilon$
respectively),
$$\phi_{gen}(g^*\E)\le
\phi_{gen}(\alpha^*g^*\E) \le \phi((\delta\circ\kappa)^*\E)
\le \phi(\E) + d.$$
\end{proof}

\section{Asymptotic regularity}\label{sec:asymptotic-regularity}

Fix a very  ample line bundle $O_X(1)$ on a projective variety $X$. 
Recall that a coherent sheaf $\F$ on $X$ is $m$-regular
\cite{mumford} provided that $H^i(\F(m-i)) = 0$ for $i>0$.
The regularity $reg(\F)$ of a sheaf $\F$ is the
least $m$ such that $\F$ is $m$-regular.
Let $Reg(X)=max(1,reg(O_X))$. Although we will not need this,
it is worth remarking that $Reg(X) = reg(O_X)$
unless $(X, O_X(1))$ is a projective space with an ample line bundle of
degree $1$.

\begin{lemma}
Let $\F$ be  $0$-regular coherent sheaf, then it is globally generated and 
$$ ker[H^0(\F)\otimes O_X\to \F]$$
is $Reg(X)$-regular.
\end{lemma}

\begin{proof}
The global generation of $0$-regular sheaves is due
to \cite[p. 100]{mumford}.
Let $R= Reg(X)$ and $K=  ker[H^0(\F)\otimes O_X\to \F]$.
By definition $R \ge 1$. We have
an exact sequence
$$0\to K(R-i)\to H^0(\F)\otimes O_X(R-i)\to \F(R-i)\to 0.$$
From the long exact of cohomology groups and $R+1$-regularity of
$O_X$ and the $R$-regularity
of $\F$ [loc. cit.], we can conclude that
$H^i(K(R-i))= 0$ for $i > 1$.
As the multiplications
$$H^0(\F)\otimes H^0(O_X(i))\to H^0(\F(i))$$
are surjective  for $i \ge 0$ [loc. cit.],
$H^1(K(R-1))$ injects into $H^1(O_X(R-1)) = 0$.
Therefore $K$ is $R$-regular.
\end{proof}

\begin{cor}\label{cor:syz}
Let $\F$ be $m$-regular, then for any $N \ge 0$, there
exist vector spaces $V_i$ and a resolution
$$V_N\otimes O_X(-m-NR)\to\ldots V_1\otimes O_X(-m-R) \to
V_0\otimes O_X(-m)\to \F\to 0$$
where $R= Reg(X)$.
\end{cor}

\begin{proof}
After replacing $\F$ by $\F(m)$, we may assume that $m=0$.
Therefore $\F$ is generated by its global sections
 $V_0 = H^0(\F)$. Let $K_0= \F(-R)$ and
$K_1$ be the kernel of the surjection
$V_0\otimes O_X\to K_0(R)$.
Then $K_1(R)$ is $0$-regular, so we can continue the above
process indefinitely and define vectors spaces $V_i$ and sheaves $K_i$
which fit into exact sequences
$$0\to K_{i+1}\to V_i\otimes O_X\to K_i(R)\to 0.$$
After tensoring these with $O_X(-iR)$, we can splice these
sequences together to obtain the desired resolution.
\end{proof}

\begin{lemma}\label{lemma:cpestim1} Let $\E$ be a coherent sheaf
on $X$, and
let $n$ be the greatest integer strictly less than $-reg(\E)/Reg(X)$.
Then 
$$\phi(\E) \le max(dim X-n-1, 0).$$
In particular, $\E$ is \cp if $reg(\E) < -Reg(X)(dim X -1)$.
\end{lemma}

\begin{proof}
Let $m = reg(\E)$, $R= Reg(X)$ and $d = dim X$. We may assume $d-n-1 > 0$, 
otherwise the lemma is trivially true.
By corollary \ref{cor:syz}, there exists a resolution
$$0\to \E_{n+1}\to \E_n\to \ldots \E_0\to \E\to 0$$
where $\E_i = V_i\otimes O_X(-m-iR)$ for $i \le n$,
and 
$$\E_{n+1} = ker[V_n\otimes O_X(-m-nR)\to V_{n-1}\otimes O_X(-m-(n-1)R)].$$
When $i \le n$, we have $-m-iR> 0$, therefore  $\phi(\E_i) = 0 \le i + d -n-1$
by lemma \ref{lemma:ampleline}.
Also $\phi(\E_{n+1}) \le d = (n+1) + d-n-1$. 
Consequently the lemma follows from theorem \ref{thm:first} (2).
\end{proof}

Suppose that $char\, k = p > 0$. Let 
$$ minreg(\E) = {\inf}_n \,\{reg\,\E^{(p^n)}\}.$$

\begin{cor}\label{cor:cpestim}
For any coherent sheaf $\E$,
$$\phi(\E) \le max(dim X-n-1, 0),$$
where $n$ is the greatest integer strictly less than $-minreg(\E)/Reg(X)$.
\end{cor}

\begin{proof} Apply the lemma to all powers
$\E^{(p^n)}$.
 \end{proof}

When $char\, k = p > 0$, we
define  the  {\em asymptotic regularity }
$$ areg(\E) = {\lim \sup}_n \, reg\, \E^{(p^n)}.$$
Of course $minreg(\E)\le areg(\E)$, but equality will usually fail.
For example, $minreg(O_X(-1)) < areg(O_X(-1)) = \infty$.
When $char\, k = 0$, define $areg(\E)$ to be the infimum of
$\sup_q[areg(\E_q)]$ over all thickenings 
of $(X,\E, O_X(1))$.
In other words, $areg(\E) \le  m$ if and only if 
$areg (\E_q) \le m$ for almost all $q$ for a given thickening.

\begin{lemma}\label{lemma:cpcrit}
($char\, k = p$)
Let $\E$ be a coherent sheaf. The following
statements are equivalent
\begin{enumerate}

\item $\E$ is $F$-ample.

\item $areg(\E) = -\infty$.

\item  $minreg(\E) < Reg(X)(dim X -1)$

\end{enumerate}
\end{lemma}

\begin{proof}
If $\E$ is $F$-ample, then clearly $reg(\E^{(p^n)})\to -\infty$ which
is the content of 2. The implication $2\Rightarrow 3$ 
follows from the inequality $minreg(\E)\le areg(\E)$.
 The implication $3\Rightarrow 1$ 
follows from corollary \ref{cor:cpestim}.
\end{proof}

\begin{cor} Conditions (1) and (2) are equivalent
 in characteristic $0$.
\end{cor} 

In any characteristic,
call $\E$ {\em F-semipositive }(with respect to $O_X(1)$)
if and only if $areg(\E) < \infty$.
We will see, shortly, that this notion is independent of the
choice of $O_X(1)$. The previous lemma shows that an
$F$-ample sheaf is $F$-semipositive.

\begin{lemma}\label{lemma:regularityestimate}
 Let $N+1\ge d=\dim X$. If 
$$\E_N\to \E_{N-1}\to\ldots \E_0\to \E\to 0$$
is an exact sequence of coherent sheaves on $X$,
then 
$$reg(\E)\le \max\{reg(\E_0),reg(\E_1)-1,\ldots reg(\E_{d-1})-d\}$$
\end{lemma}

\begin{proof} Extend this to a sequence
$$0\to \E_{N+1}\to \E_N\to\ldots\E_0\to \E\to 0.$$
The regularity estimate follows from lemma~\ref{lemma:one}
and the fact that $m$-regular sheaves are $m'$-regular
for all $m'\ge m$ \cite[p. 100]{mumford}.
\end{proof}

\begin{prop}\label{prop:semipos} 
Let $f:X\to Y$ be a morphism of  projective varieties.
Assume that $Y$ is equipped with a very ample line bundle
$O_Y(1)$. Let $\E$ be a coherent sheaf on $Y$ which is
  F-semipositive  with  respect  to $O_Y(1)$. If
${\mathcal T}or_i^{f^{-1}O_Y}(O_X, \E) = 0$ for all $i > 0$ (e. g.
if  $f$ is flat, or $\E$ is locally free), then
$f^*\E$ is F-semipositive with respect to $O_X(1)$.
\end{prop}

\begin{proof}
We give the proof in positive characteristic.
By hypothesis and  corollary \ref{cor:syz}, there exists a resolution.

$$V_N\otimes O_Y(-m-NR)\to\ldots V_1\otimes O_Y(-m-R) \to
V_0\otimes O_Y(-m)\to \E^{(p^n)}\to 0$$
where the constants $m, R, N>>0$ can be chosen
independently of $n$. This stays exact after applying $f^*$
by our assumptions. Therefore the regularity
of $f^*(\E^{(p^n)}) = (f^*\E)^{(p^n)}$ stays bounded as
$n\to \infty$ by lemma~\ref{lemma:regularityestimate}.
\end{proof}

\begin{cor}
Let $O_X(1)'$ be another very ample line bundle on $X$,
then a sheaf $\E$ is F-semipositive with respect to 
$O_X(1)$ if and only if it is F-semipositive with respect to 
$O_X(1)'$.
\end{cor}

\begin{proof}
Apply the proposition to the identity map.
\end{proof}

Recall that a locally free sheaf $\E$ on $X$
 is nef (or numerically semipositive)
if  for any curve $f:C\to X$, any quotient of $f^*\E$
has nonnegative degree. In characteristic $0$, it is convenient
to introduce an ostensibly stronger property: $\E$ is arithmetically
nef if there is a thickening $(\tilde X,\tilde \E)$ over $Spec\, A$
such that the restriction  of $\tilde \E$ to the fibers are nef.
To simplify the statements, we define arithmetically
nef to be  synonymous with nef in positive characteristic.
Further discussion of these matters can be found in the appendix.
The name $F$-semipositive stems from the following:

\begin{lemma}\label{lemma:Fsemiposisnef} 
If $\E$ is an $F$-semipositive locally free sheaf, then
it is arithmetically nef.
\end{lemma}

\begin{proof}
By definition, we may work over a field of characteristic $p>0$.
Suppose that $\F$ is a quotient of $f^*\E$ with negative
degree. This implies that $deg(\F^{(p^n)})\to -\infty$
as $n\to \infty$.
By proposition \ref{prop:semipos}, $f^*\E$ is F-semipositive,
which implies that there is a fixed line bundle $L$
such that $\E^{(p^n)}\otimes L$ is globally generated for
all $n$. Therefore $\F^{(p^n)}\otimes L$ is globally
generated for all $n$ which implies that $deg(\F^{(p^n)})$
is bounded below. This is a contradiction.
\end{proof}

For line bundles, the  converse is given by
proposition~\ref{prop:nefline}. However, it fails for higher rank
(example~\ref{ex:amplenotFample}).

\section{Tensor products}

\begin{thm}\label{thm:tensor1}
Let $\E$ and $\F$ be two vector bundles on a smooth projective variety
$X$, then
$$\phi(\E\otimes \F) \le \phi(\E)+\phi(\F).$$
\end{thm}

\begin{proof}
Assume that $k$ is a field of characteristic $p>0$.
  Let $Y = X\times X$ and let $p_i:Y\to X$ denote the projections.
Given two coherent sheaves $\E_i$ on $X$, let 
$\E_1\boxtimes \E_2 = p_1^*\E_1\otimes p_2^*\E_2$.
Choose a very ample line bundle $O_X(1)$ on $X$, then
$L = O(1)\boxtimes O(1)$ is again very ample. Let $\Delta\subset X$ be
the diagonal. Choose $\nu >>0$. By corollary \ref{cor:syz},
we can construct a resolution 
\begin{equation}\label{eq:Diagres}
0\to \G_{\nu+1}\to \G_\nu\to \ldots \G_0\to O_\Delta\to 0
\end{equation}
where $\G_i = V_i\otimes L^{\otimes a_i}$ for $i\le \nu$.

The Frobenius map $F_Y = F_X\times F_X$. Thus using K\"unneth's
formula \cite[III, 6.7.8]{ega}, for any $b$ we get
\begin{eqnarray*}
H^i(\G_j\otimes L^{b}\otimes F_Y^{N*}(\E\boxtimes \F))
&=& V_j\otimes H^i(\E^{(p^N)}(b+a_j)\boxtimes \F^{(p^N)}(b+a_j))\\
&=&  V_j\otimes \bigoplus_{c+d=i} \, H^c(\E^{(p^N)}(b+a_j))\otimes
H^d(\F^{(p^N)}(b+a_j))\\
&=& 0
\end{eqnarray*}
for $i>\phi(\E) + \phi(\F)$, $j\le \nu$ and $N>>0$. 
Tensoring (\ref{eq:Diagres})
by $L^{b}\otimes F_Y^{N*}(\E\boxtimes \F)$ and applying 
lemma \ref{lemma:one} shows that 
\begin{eqnarray*}
H^i((\E\otimes \F)^{(p^N)}(2b)) &=&
H^i(O_\Delta\otimes L^{b}\otimes F_Y^{N*}(\E\boxtimes \F))\\
& =& 0
\end{eqnarray*}
for  $i>\phi(\E) + \phi(\F)$ and $N>>0$.
Thus corollary \ref{cor:FamplViaLineb} gives the desired bound on
$\phi(\E\otimes \F)$.

If $char\, k = 0$, then we can carry out the above argument
on the fiber of some thickening.
\end{proof}

\begin{cor}\label{cor:tensorample}
  The tensor product of two $F$-ample bundles  is $F$-ample.
\end{cor}

\begin{cor}\label{cor:logtensor}
Let $D$ and $E$ be  reduced effective divisors
such that $D+E$ has normal crossings. Suppose that 
 $\E$ and $\F$ are a pair of vector bundles with critical
divisors $\Delta$ and $\Xi$ along $D$ and
$E$ respectively. If $\Delta+\Xi$ is strictly fractional, i.e. has
all its multiplicities less then $1$, then
$$\phi(\E\otimes \F, D+E)\le \phi(\E,D) +\phi(\F,E).$$
\end{cor}

\begin{proof}
We can find $n,m$ such that
  $\phi(\E,D) = \phi(\E^{(p^n)}(-D'))$ and  $\phi(\F,E) =
  \phi(\F^{(p^m)}(-E'))$ where $D' = p^n\Delta$
and $E'= p^m\Xi$. After replacing $E'$ by $p^{n-m}E'$,
or the other way around, we can assume that $m=n$. Therefore
$$(\E\otimes \F)^{(p^n)}(-D'-E') = \E^{(p^n)}(-D')\otimes
\F^{(p^n)}(-E')
$$ has $F$-amplitude bounded by the sum.
\end{proof}

\begin{remark} If $D$ and $E$ are disjoint, then the conditions
on the critical divisors are automatic.
\end{remark}


  


\begin{thm}\label{thm:tensor} 
Let $\E$ and $\F$ be two  coherent sheaves on $X$ such that one of them is
locally free and
 $\E$ is F-semipositive, then

$$\phi(\E\otimes \F) \le \phi(\F).$$

\end{thm}

\begin{proof} Assume that $char\, k = p >0$.
 Let $m=areg(\E)$, $N>>0$ and $R= Reg(X)$.
Then $reg(\E^{(p^\mu)}) \le m$ for all but finitely many
$\mu$.  Given a locally free sheaf ${\mathcal G}$,
choose $\mu_{0}$, so that 
$$H^i(\F^{(p^\mu)}\otimes {\mathcal G}(-m-jR))=0$$
for all $\mu> \mu_{0}$, $i> \phi(\F)$, and $j = 0,1,\ldots N$.
By increasing $\mu_0$, if necessary, we can assume that
 $\E^{(p^\mu)}$ is $m$-regular when $\mu> \mu_o$.
From corollary \ref{cor:syz}, we
 obtain a resolution
$$0\to \E_{N+1}\to \E_N\to \ldots \E_0\to \E\to 0$$
where $\E_i= V_i\otimes O_X(-m-iR)$ for $i \le N$. Tensoring this by $\F^{(p^\mu)}\otimes {\mathcal G}$
and applying lemma \ref{lemma:one} shows that
$$H^i((\E\otimes\F)^{(p^\mu)}\otimes {\mathcal G}) = 0$$
when $\mu> \mu_{0}$ and $i > \phi(\F)$.

If $char\, k = 0$, then we can carry out the above argument
on the fiber of some thickening.
\end{proof}

We can refine  corollary \ref{cor:tensorample}.

\begin{cor}\label{cor:tensorcp}
 The tensor product of  an \cp vector bundle and an F-semipositive 
vector bundle is \cp.
\end{cor}

\section{Characterization of $F$-ample sheaves  on special varieties}

It is possible to give an elementary characterization of
$F$-ampleness for curves and projective spaces. Recall that a vector
bundle $\E$ over a variety $X$ defined over a field $k$ of characteristic
$p$  is $p$-ample \cite{hartshorne} if for any coherent sheaf $\F$
there exists  $n_0$ such that $\E^{(p^n)}\otimes \F$ is globally
generated for all $n\ge n_0$.

\begin{lemma}\label{lemma:Fample2pample} 
    An $F$-ample vector bundle $\E$ is $p$-ample.
\end{lemma}

\begin{proof} 
 Suppose $\F$ is a coherent sheaf. Since the regularity
 of the sheaves $\E^{(p^n)}\otimes \F\to -\infty$, these
 sheaves are globally generated for $n>>0$.
 \end{proof}

\begin{cor}  An $F$-ample vector bundle $\E$ is ample. \end{cor}
    
\begin{proof} \cite[6.3]{hartshorne} \end{proof}

As we will see the converse to both statements fail in
general. 

\begin{lemma}\label{lemma:phipample}
    If $\E$ is a $p$-ample vector bundle on a projective
variety $X$, then $\phi(\E) < dim\, X$.
\end{lemma}

\begin{proof} 
 Choose a coherent sheaf $\F$.   
 Let $L$ to be the $N$th power of an ample line bundle,
 chosen large  enough so that $H^i(\F\otimes L) = 0$ for $i>0$. Then for
 all $n>>0$, $\E_n=\E^{(p^n)}\otimes L^{-1}$ is globally generated. Therefore
 $$H^0(\E_n)\otimes \F\otimes L\to \E^{(p^n)}\otimes \F$$
 is surjective. It follows that the top degree
 cohomology of the right hand vanishes for $n>>0$.
\end{proof}

This leads to a complete characterization for curves.

\begin{prop}\label{prop:curveFample-ample}
Let $\E$ be a coherent sheaf over a smooth projective curve $X$
defined over a field $k$,
then  the following are equivalent
\begin{enumerate}
    \item $\E$ is $F$-ample.
    \item $\E/torsion$ is $p$-ample when $char\, k= p$.
    \item $\E/torsion$ is ample.
 \end{enumerate}
\end{prop}

\begin{proof}
  Since $\E$ is a direct sum of $\E/torsion$ with the torsion
  part, we may assume that $\E$ is a vector bundle.
 Suppose $char\, k=p$. Then the equivalence of the  first
 two statements follows from lemmas \ref{lemma:Fample2pample} 
and \ref{lemma:phipample}. The equivalence of the last
two from \cite[6.3, 7.3]{hartshorne}. 

If $char\, k = 0$, we can deduce the equivalence
of (1) and (3) from the previous cases, because ampleness
is an open condition.      
\end{proof}

Now, we turn to projective space.
For integers $a\le b$, let
 $[a,b] = \{a, a+1, \ldots , b\}$.

\begin{thm}\label{thm:proj}
Let $\E$ be a coherent sheaf on the projective space
$\PP^n_k$, then
$$\phi(\E) = min\, \{i_0\,|\, H^i(\PP_k^n, \E(j)) = 0,\,
\forall i> i_0, \forall j\in [-n-1, 0]\}$$
In particular, $\E$ is $F$-ample
if and only if $H^i(\E(j)) = 0$ for all $j\in [-n-1, 0]$
and $i > 0$.
\end{thm}

This leads to a characterization of $F$-ample bundles
on  $\PP^n= \PP_k^n$. For a slightly different characterization,
see \cite[sect. 4]{mi}.
The key step in the proof of theorem \ref{thm:proj}
is the following proposition.

\begin{prop}\label{prop:split}
 Let $\pi:\PP^n\to \PP^n$ be a finite morphism,
and let $d$ be the degree of $\pi^*O(1)$.

\begin{enumerate}

\item For each $i$,
$\pi_* O(i)$ is a direct sum of line bundles.

\item If $-d -n -1 < i < d$ then  $\pi_* O(i)$
is a sum of line bundles of the form $O(l)$ with $l\in [-n-1, 0]$.

\item There exists a constant $C$ depending only on $n$, such that
if $d > C$, then for each $O(l)$ with
$l\in [-n-1, 0]$, $O(l)$ occurs as  a component  of  $\pi_*O(i)$
for some $i\in [-n-1, 0]$.
\end{enumerate}

\end{prop}

\begin{proof} 
Since $\pi$ is finite, 
$$H^a((\pi_*O(i))\otimes O(j))= H^a(O(i+ dj))$$
for all $a,i$ and $j$. In particular, these groups vanish for
all $0< a < n$ and all $j$. Therefore $\pi_*O(i)$ splits
into a sum of line bundles by a theorem of Horrocks
\cite{horrocks}.

If $-d -n -1 < i < d$, then 
$$H^n((\pi_*O(i))\otimes O(1))= H^n(O(i+ d))= 0$$
and
$$H^0((\pi_*O(i))\otimes O(-1))= H^0(O(i- d))= 0.$$
This implies the second statement.

Let $p_m(x) = \left(
\begin{array}{c}x+m\\ m\\ \end{array} \right)$
and $(\Delta_x p)(x,y,\ldots) = p(x,y,\ldots) -p(x-1,y,\ldots)$.
Note that $\Delta_x p_m(x) = p_{m-1}(x)$.

Choose $-n-1 \le i \le 0$.
 Let us write
$$\pi_*O(i) = \bigoplus_{l} O(l)^{\oplus f(l,i)}.$$
By comparing cohomology of $\pi_*O(i)$ and $O(i)$,
we see that
$S=\{l\,|\, f(l,i)\not= 0\}$ is contained in $[-n, 0]$
if $i=0$, $S$ is  contained in
$[-n, -1]$ if $-n-1 < i < 0$, and 
$S$ is contained in $[-n-1,-1]$ if $i = -n-1$.
Furthermore $f(0,0) = f(-n-1, -n-1)=1$, and this
shows that the proposition holds true for $l= 0,-n-1$.
We now assume  that $-n\le l \le 0$.
Tensoring $\pi_*O(i)$ by $O(x)$ and computing
Euler characteristics  yields:
\begin{equation}\label{eq:fp}
\sum_l f(l,i)p_n(x+l) = p_n(dx+i).
\end{equation}
Setting $x=0$ yields
$$f(0,i) = p_n(i).$$
Applying $\Delta_x$ to equation (\ref{eq:fp}) and setting $x=0$ yields
$$(-1)^{n-1}f(-n,i) = (p_n(i) - f(0,i)) - p_n(i-d)$$
hence
$$ f(-n,i)=(-1)^n p_n(i-d).$$
Applying $\Delta_x^2$ to equation (\ref{eq:fp}) and setting
$x=0$ yields
$$(-1)^{n-2}(n-1)f(-n, i) + (-1)^{n-2}f(-n+1, i) 
= (p_n(i) - f(0,i)) - 2p_n(i-d) + p_n(i-2d)$$
hence 
$$f(-n+1, i) = (-1)^{n-1}(n+1)p_n(i-d) + (-1)^{n-2}p_n(i-2d).$$
Doing this repeatedly gives a formula for each
$f(l,i)$, with $-n-1< l <  0$, 
which is a nonzero polynomial of degree at
most $n$ in $i$ and $d$. Choosing a specific $d >> 0$ (for $n$ fixed)
forces $f(l,i)$ to be a nonzero polynomial in $i$ of degree
at most $n$. Therefore $f(l,i)\not= 0$ for some $i$
in the range $-n\le i \le  0$.
\end{proof}

We need the following (presumably well known) version
of the projection formula.

\begin{lemma}\label{lemma:projform}
 If $\pi:X\to Y$ is a finite map of quasiprojective
schemes, then $\pi_*(\pi^*\E\otimes \F) \cong \E\otimes \pi_*\F$.
\end{lemma}

\begin{proof} Choose a resolution 
$\E_1\to \E_0\to \E\to 0$ by vector bundles $\E_i$.  There
is a diagram
$$
\begin{array}{ccccccc}
\E_1\otimes \pi_*\F& \to& \E_0\otimes \pi_*\F& \to& \E\otimes \pi_*\F&
\to& 0\\
\downarrow&&\downarrow&&\downarrow&&\\
\pi_*(\pi^*\E_1\otimes \F)& \to& \pi_*(\pi^*\E_0\otimes \F)& \to&
\pi_*(\pi^*\E\otimes \F)& \to& 0\\
\end{array}
$$
where the first two vertical arrows are isomorphisms by the usual
projection formula. This implies that the third arrow is
also an isomorphism.
\end{proof}

\begin{proof}[Proof of theorem \ref{thm:proj}]
As usual, we prove the result in positive characteristic;
the characteristic zero case is a formal consequence.
To begin with, we show $H^i(\E(j)) = 0$ for $i>\phi(\E)$ and
$j\in [-n-1, 0]$. First assume $char\, k = p> 0$.
Choose $m>>0$, then $O(j)$ is a
direct summand of some $F^m_*O(l)$ for $l\in [-n-1, 0]$
by the previous proposition.
Therefore by the projection formula (lemma \ref{lemma:projform}),
$$H^i(\E(j)) \subseteq H^i(\E\otimes F^m_*O(l))
= H^i(\E^{(p^m)}\otimes O(l)) = 0.$$

Conversely, suppose that $H^i(\E(j)) = 0$ for
all $i > i_0$ and $j\in [-n-1, 0]$. 
For any
integer $l$, we can choose $m >> 0$ so that
$F^m_*O(l)$ is a direct of line bundles 
$O(j)$ with $j \in [-n-1, 0]$. Therefore
$$ H^i(\E^{(p^m)}\otimes O(l))=H^i(\E\otimes F^m_*O(l)) = 0$$
for $i > i_0$. Since any coherent sheaf $\F$ on $\PP^n$
has a finite resolution by direct sums of line bundles, this
shows that $H^i(\E^{(p^m)}\otimes \F) = 0$ for $m>>0$ and $i>i_0$
(by the same argument as in the proof of
corollary~\ref{cor:FamplViaLineb}).
\end{proof}

The theorem yields improvements on the
regularity estimates of the previous section.
\begin{cor} An $F$-ample sheaf $\E$ on $\PP^n$
is  $(-1)$-regular; in particular $\E(-1)$ is globally generated.
\end{cor}

\begin{ex}\label{ex:amplenotFample}
 The tangent bundle $T$ of $\PP^n$ is ample and in fact
  $p$-ample in positive characteristic,
 but $T$ is not $F$-ample if $n\ge 2$, because
$H^{n-1}(T(-n-1)) = H^1(\Omega^1)^*\not=0$.
The bundle $T(-1)$ is globally generated, and therefore nef.
However, it cannot be $F$-semipositive, since otherwise $T$
would be $F$-ample by corollary~\ref{cor:tensorcp}.
\end{ex}

\begin{ex}\label{ex:restr} Let $X$ be a projective variety with an ample
line bundle $L$. Embed $i:X\hookrightarrow \PP^n$ using a large multiple
of $L$. Then $L$ is \cp but $i_*L$ is not, because the conclusion
of the above corollary fails. Therefore theorem~\ref{thm:first} (5)
fails for non\'etale finite maps.
\end{ex}

\begin{cor} A vector bundle on $\PP^2$ is
$F$-ample if and only if it
is isomorphic to a sum of the
form $E\oplus O(1)^{\oplus N}$ where $E$ is $(-2)$-regular. 
\end{cor}

\begin{proof}
A  direct sum of a $(-2)$-regular sheaf and 
a bunch of $O(1)$'s satisfies the  conditions
of the theorem by \cite[p. 100]{mumford}.

Now suppose that $V$ is a \cp vector bundle.
It is $(-1)$-regular by the previous corollary,
and therefore $V(-1)$ is generated by global sections.
Suppose that $H^2(V(-4))\not= 0$. Then by Serre duality,
there is a nonzero morphism $V(-1)\to O$ and let
$V'$ be the kernel twisted by $O(1)$. Since the map 
$H^0(V(-1))\otimes O\to O$
must split, it follows that  the map $V(-1)\to O$ also splits.
Therefore $V'$ is again $(-1)$-regular, so we can continue splitting off
copies of $O(1)$ from $V$ until we arrive at a summand $E$ with
$H^2(E(-4)) = 0$. Since we also have $H^1(E(-3)) = 0$,
 it follows that $E$ is $(-2)$-regular.
\end{proof}

\section{$F$-amplitude of ample bundles}

As we have seen, ample vector bundles need not be $F$-ample.
However, we do have an estimate on their amplitude, at least in 
characteristic $0$.

\begin{thm}\label{thm:keyestimate}
 Let $X$ be a projective variety over
a field of characteristic $0$ and let $\E$ be an ample vector bundle
of rank $r$ on $X$.  Then $\phi(\E)< r$. 
\end{thm}

Keeler has found that the inequality $\phi(\E)< dim\, X$ also 
holds for ample vector bundles (proposition~\ref{prop:ample-dimension-bound}).

Before giving the proof, we need to review some results
from (modular) representation theory. We will choose
our notations consistent with those of 
\cite{arapura}.
 Let $A$ be a
commutative ring and $E=A^r$.
Fix a partition $\lambda = (\lambda_1\ge\lambda_2\ge\ldots)$
of weight $|\lambda| = \sum_i\lambda_i$. The Schur power
$\BS^\lambda(E)$ can be constructed as the space of global
sections of a line bundle associated to $\lambda$
over the scheme $Flag(E)$ of flags on $E$. 
A more elementary construction is possible;
$\BS^\lambda(E)$ can be defined as 
a quotient of $\E^{\otimes|\lambda|}$ by an explicit
set of relations involving $\lambda$ \cite[8.1]{fulton}.
This quotient map can be split using a Young symmetrizer
when $A$ contains $\Q$, but not in general.
The paper \cite{carter-lustig} gives essentially
a dual construction; their Weyl module $E_\lambda$ coincides
with our  $\BS^{\lambda'}(E)^*$
where $\lambda'$ is the conjugate partition (see \cite[p 251]{jantzen}
for the comparison with the first construction). 
 We will need to make the initial description more explicit.
Let $\pi_k:Flag(E)\to Grass_k(E)$ be the canonical map to
the Grassmannian of $k$-dimensional {\em quotients} of $E$, and $i_k$ its
Pl\"ucker embedding. Then
$$\BS^\lambda(E) = 
H^0(Flag(E), L_\lambda)$$
where $a_i = \lambda_i-\lambda_{i+1}$ and
$$L_\lambda =\bigotimes_k\, \pi_k^*i_k^*O_{\PP(\wedge^kE)}(a_k).$$

At the two extremes, $\BS^{(n)}(E) = S^n(E)$ and
$\BS^{(1,1,\ldots 1)}(E) = \wedge^i(E)$ where $i$ is the length
of the string. These Schur powers turn out to be free
$A$-modules, and since the constructions are functorial, they
carry  $GL_r(A)$ actions.
When $A=k$ is a field of characteristic $0$,
 the $GL_r(k)$-modules $\BS^\lambda(E)$ are all irreducible.
This is no longer true when $A=k$ is a field of characteristic
$p>0$. For example, the symmetric power $S^p(E)$ contains
a nontrivial submodule $E^{(p)}$ which is the representation
associated to the $p$th power map $GL_r(k)\to GL_r(k)$.
This inclusion can be extended to a resolution:

\begin{thm}(Carter-Lusztig  \cite[pg 235]{carter-lustig})
If $k$ is a field of characteristic $p>0$. Then 
there exists an exact sequence of $GL_r(k)$-modules
$$0\to E^{(p)}\to \BS^{(p)}(E)\to \BS^{(p-1, 1)}(E)\to
 \BS^{(p-2, 1, 1)}(E)\to\ldots \BS^\lambda(E)\to 0$$
where  $\lambda=(p-\min(p-1,r-1), 1,1\ldots)$.
\end{thm}

These constructions are easy to globalize to the case
where $E$ is replaced by a vector bundle $\E$.
In this case, the two meanings of $\E^{(p)}$ agree.

\begin{cor}
If $\E$ is a vector bundle over a scheme $X$ defined
over a field $k$ of characteristic $p>0$, then 
there exists a resolution of $\E^{(p)}$ as above.
\end{cor}

Suppose we are in the situation of theorem \ref{thm:keyestimate}.
Choose a line bundle $M$ on $X$.
Let 
$(\tilde X, \tilde \E, \tilde M)$
 be an arithmetic thickening over $ Spec\, A$.
After shrinking $Spec\, A$ we can assume that $\E$ is locally free.
Then we have vector bundles $\BS^\lambda(\tilde \E)$ 
over $\tilde X$. Fix a partition $\lambda$.
Then for any natural number
$N$, we get a new partition $(N)+\lambda = (N+\lambda_1, \lambda_2,\ldots)$.

\begin{lemma}\label{lemma:HiXqBS} 
With notation and assumptions as above
(specifically that $\E$ is ample),
there exists  an integer $N_0$ such that 
$$H^i(\tilde X, \BS^{(N)+\lambda}(\tilde \E)\otimes \tilde M) = 0$$
for $N\ge N_0$.
\end{lemma}

\begin{proof}
Let $\pi:Flag(\E)\to X$ be the bundle of flags on $\E$.
To simplify notation, we will write $\tilde M$ instead of
$\pi^*\tilde M$.
The fibers of $\pi$ are partial flag varieties. The higher
cohomology groups of $L_{(N)+\lambda}$  along these fibers are zero by 
Kempf's vanishing theorem (see for example \cite[II, 4.5]{jantzen}).
 Therefore the higher direct images vanish, and consequently
the Leray spectral sequence yields isomorphisms
$$H^i(\tilde X, \BS^{(N)+\lambda}(\tilde \E)\otimes \tilde M) \cong
H^i(Flag(\E), L_{(N)+\lambda}\otimes \tilde M).$$
Let $\pi_1:Flag(\E)\to \PP(\E)$ be the canonical projection.
For reasons similar to those above, there are isomorphisms
$$H^i(Flag(\E), L_{(N)+\lambda}\otimes \tilde M)
 \cong  H^i(\PP(\E), \pi_{1*}L_{(N)+\lambda}\otimes\tilde M).$$
By the projection formula, the right hand side is the cohomology
of $O_{\PP(\E)}(N)\otimes \pi_*(L_\lambda)\otimes \tilde M$. Since $O(1)$
is ample, these groups vanish for $N>>0$ and $i>0$.
\end{proof}

\begin{cor}\label{cor:HiXqBS}
With the notation of section 1, given $a > 0$, there exists
$N_0$
$$H^i( X_q, \BS^{(N)+\lambda}(\E_q)\otimes M_q^{\otimes n}) = 0$$
for all $i>0$ , $0\le n\le a$,
$N\ge N_0$ and closed points $q\in Spec\, A$.
\end{cor}

\begin{proof}
\cite[III. 12.9]{hartshorne2}
\end{proof}

\begin{proof}[Proof of theorem \ref{thm:keyestimate}]
Choose $M = O_X(-1)$ with $O_X(1)$ very ample.
Let $C << 0$ be a constant. 
By corollary \ref{cor:HiXqBS}, there exists
a $N_0$ such  that the sheaves
$ \BS^{(N-i, 1,\dots 1)}(\E_q))$ ($0\le i< r$) have regularity less
than $C$ for all $N\ge N_0$ and all closed points $q\in Spec\, A$.
In particular,  there exits  a nonempty open set
$U\subset Spec\, A$ such that 
$$reg( \BS^{(p(q)-i, 1,\dots 1)}(\E_q)) < C$$
for all closed $q\in U$ ($p(q) = char\, A/q$).
By lemma \ref{lemma:cpestim1}, these sheaves are $F$-ample.
Then the Carter-Lusztig resolution (which has
length bounded by $r=rank(\E)$) together with lemma
\ref{lemma:one} shows that $\phi(\E)< r$.
\end{proof} 

\begin{cor}\label{cor:keyestimate}
Let $\E_i$ be ample vector bundles. Then $$\phi(\E_1\otimes\E_2\otimes\ldots
\E_m) < rank(\E_1)+rank(\E_2)+\ldots rank(\E_m).$$
\end{cor}

The analogue for a pair is the following:

\begin{thm}\label{thm:logkeyestimate}
 Let  $X$ be a smooth projective variety defined
over a field of characteristic $0$, and let $\E$ be a vector bundle
of rank $r$ on $X$. Suppose there exists a reduced normal crossing
divisor $D$,  a positive integer $n$, and
a divisor $0\le D'<_{strict} nD$ such that $S^n(\E)(-D')$ is ample.
Then $\phi(\E, D) < r$.
\end{thm}

\begin{remark} The hypothesis amounts to the condition that
    the ``vector bundle'' $\E(-\Delta)$ is ample for some
fractional effective $\Q$-divisor $\Delta =\frac{1}{n}D'$
supported on $D$.
\end{remark}

\begin{proof}
The proof is very similar to the previous one, so we will just
summarize the main points. Let
$M = O_X(-1)$ with $O_X(1)$ very ample. Choose a thickening of $(X,\E,D,M)$.
A small modification of corollary \ref{cor:HiXqBS} shows
that the regularity of the sheaves 
$$ \BS^{(Nn+j-i, 1,\dots 1)}(\E_q)(-ND'_q),\> 0\le i< r,\, 0<j<n$$ 
can be made less than a given $C$ for all $N$ greater than some
$N_0$ depending on $C$. All but finitely primes are of the form
$Nn+j$ for $N$ and $j$ as above. Thus the above sheaves will be $F$-ample
for  almost all $q$. The Carter-Lusztig resolution shows that
$\phi(\E^{(p(q))}(-ND')) < r$ which implies the theorem.
\end{proof}





\section{An $F$-ampleness criterion}

The notion of geometric positivity was introduced in \cite{arapura}.
Although the methods are very different, there appear to be some parallels
between $F$-ampleness and geometric positivity.
The following result is an analogue of  [loc. cit., cor. 3.10].

\begin{thm}\label{thm:strongsemistable}
  Let $\E$ be a rank $r$ vector bundle on a smooth projective variety
$X$ such that $det(\E)$ is ample and
$S^{rN}(\E)\otimes det(\E)^{-N}$ is globally generated for some
$N>0$ prime to $char\, k$. Then $\E$ is $F$-ample.
\end{thm}

\begin{remark} The hypothesis that 
$S^{rN}(\E)\otimes det(\E)^{-N}$ is globally generated implies that
$\E$ is strongly semistable in the sense of \cite[p. 247]{arapura}.
We leave it as an open problem to determine whether $F$-ampleness
follows only assuming  strong semistability of $\E$ and ampleness
of $det(\E)$.
\end{remark}

As usual all the work will be in characteristic $p>0$. Let $q= p^n$
for some $n>0$. In this section we will modify our previous conventions,
and write $F_X:X\to X$ for the absolute  $q$th power  Frobenius.
Let $\PP= \PP(\E)$ and $\PP' = \PP(\E^{(q)})$ with canonical projections
denoted by $\pi$ and $\pi'$.
Consider the commutative diagram
$$
\xymatrix{
 \PP\ar[r]_{\Phi}\ar[rd]_{\pi}\ar@/^/[rr]^{F_\PP} & \PP'\ar[r]_{\phi}\ar[d]^{\pi'} &
 \PP\ar[d]^{\pi} \\ 
  & X\ar[r]_{F_X} & X 
}
$$
where the right hand square is cartesian. $\Phi$ is the relative $q$th
power Frobenius associated to $\pi$. Let $\Sigma_N = S^{rN}(\E)\otimes det(\E)^{-N}$.

\begin{prop}\label{prop:relFsplit}
  If $\Sigma_{q-1}$ is globally generated, then 
$O_{\PP'}\to \Phi_*O_\PP$ splits.
\end{prop}

\begin{proof}
By Grothendieck duality for finite flat maps \cite[ex.III 6.10,
 7.2]{hartshorne2}, the proposition is equivalent 
to the splitting of the trace map
\begin{equation}\label{eq:Grothtrace}
tr: \Phi_*\omega_{\PP/\PP'} \to O_{\PP'}
\end{equation}

  We have
$$F^*_\PP O_\PP(1) = O_\PP(q)$$
$$\phi^* O_\PP(1) = O_{\PP'}(1),$$
therefore
$$ \Phi^*O_{\PP'}(1) = O_{\PP}(q).$$

Also
$$\omega_{\PP'/X} = O_{\PP'}(-r) \otimes(\pi')^* det(\E^{(q)}) = 
 O_{\PP'}(-r) \otimes (\pi')^*det(\E)^q$$
$$\omega_{\PP/X} = O_{\PP}(-r) \otimes \pi^* det(\E)$$
 \cite[ex. III 8.4]{hartshorne2}. Therefore
$$\Phi^*\omega_{\PP'/X}  = 
 O_{\PP}(-qr) \otimes \pi^*det(\E)^q$$
$$\omega_{\PP/\PP'}  = 
 O_{\PP}((q-1)r) \otimes \pi^*det(\E)^{1-q}.$$
Observe that
\begin{equation}\label{eq:piomegaPP}
\pi_*\omega_{\PP/\PP'} =\Sigma_{q-1}.
\end{equation}

Suppose that $0<i<r$, then using the projection formula
and the previous computations
\begin{eqnarray*}
R^i\pi'_*[(\Phi_*\omega_{\PP/\PP'})(-i)] &=&
R^i\pi_*\omega_{\PP/\PP'}(-qi)\\
&=& R^i\pi_*O_{\PP}(q(r-i)-r) \otimes \pi^*det(\E)^{1-q}\\
&=& 0
\end{eqnarray*}
Thus $\Phi_*\omega_{\PP/\PP'}$ is regular relative to $\pi'$,
and it follows that the canonical map
$$(\pi')^*\pi'_* \Phi_*\omega_{\PP'/\PP}\to \Phi_*\omega_{\PP'/\PP}$$
is surjective \cite[V, 2.2]{fulton-lang}.
By  (\ref{eq:piomegaPP}), this gives a surjection
\begin{equation}\label{eq:Sigma2omega}
(\pi')^*\Sigma_{q-1}\to \Phi_*\omega_{\PP'/\PP}
\end{equation}
Composing this with  the Grothendieck trace (\ref{eq:Grothtrace}),
gives a surjection
$$(\pi')^*\Sigma_{q-1}  \to  O_{\PP'}$$
Since $\Sigma_{q-1}$ is globally generated, there exists
a morphism $s:O_{\PP'}\to (\pi')^*\Sigma_{q-1}$ such that the composite
$O_{\PP'}\to O_{\PP'}$ is  nonzero. This corresponds
to an  element $a\in k^*=H^0(O_{\PP'})-\{0\}$.
The composite of $\frac{1}{a}s$ with (\ref{eq:Sigma2omega})
gives a splitting of (\ref{eq:Grothtrace}).
\end{proof}

\begin{cor}\label{cor:relFsplit}
 With the same  assumptions as in the proposition, $\E^{(q)}$ is a direct
summand of $S^{q}(\E)$.
\end{cor}

\begin{proof}
 By the projection formula, the canonical map 
\begin{equation}\label{eq:OPP12PhiOPPq}
O_{\PP'}(1)\to \Phi_*\Phi^*O_{\PP'}(1) = \Phi_*O_{\PP}(q)
\end{equation}
can be identified with
$$O_{\PP'}(1)\to  O_{\PP'}(1)\otimes \Phi_*O_{\PP}$$
This splits. Applying $\pi'_*$ to (\ref{eq:OPP12PhiOPPq})
yields a split injection $\E^{(q)}\to S^{q}(\E)$.
\end{proof}

\begin{proof}[Proof of theorem \ref{thm:strongsemistable}]
Choose $q=p^n\equiv 1 \> (mod\, N)$;  
$q$ can be chosen arbitrarily large.
Then $\Sigma_{q-1}$ is globally generated since
it is a quotient of $S^q(\Sigma_N)$. It follows that
$\E^{(q)}$ is a direct summand of $S^q(\E)$. Since $S^N(\E) =
\Sigma_N\otimes det(\E)^{N}$ is ample, the same is true for
$\E$. Thus $reg(\E^{(q)}) = reg(S^q(\E))\to -\infty$ as $q\to \infty$.  

This finishes the proof in characteristic $p$, the remaining case is
handled as usual.
\end{proof}

\section{The main vanishing theorem}

As  a warm up to the main theorem, we will
 extend some of the conclusions of theorem~\ref{thm:proj} 
to a more general class of spaces
called  Frobenius split  varieties \cite{mehta}.
This  means that
the map $O_{X}\to F_{*}O_{X}$ splits (actually,
we only need the ostensibly weaker property that
this map splits in the derived category).
Proposition \ref{prop:relFsplit}  implies that projective
spaces are Frobenius split. Other
examples of Frobenius split varieties include quotients
of semisimple groups by parabolic subgroups [loc. cit], and
most mod $p$ reductions of a smooth Fano variety \cite[4.11]{smith}.

\begin{prop} Suppose that
$X$ is a smooth projective variety such that
$O_X\to F_*O_X$ splits in the derived category of
$O_X$-modules. Then
$H^{i}(X,\E) = 0$ for $i > \phi(\E)$. If $\E$ is locally
free, then $H^i(X,\omega_{X}\otimes \E) = 0$
for $i > \phi(\E)$.
\end{prop}

\begin{proof}
 For the first statement, use the fact there
 is an injection
 $$H^i(X,\E)\hookrightarrow H^i(X,\E\otimes F_*O_X)$$
 because it splits by hypothesis.
 On the other hand the projection formula gives
 $$H^i(X,\E\otimes F_*O_X)\cong H^i(X,\E^{(p)})$$
 By iterating we get a  sequence of injections
 $$H^{i}(\E)\hookrightarrow H^i(\E^{(p)})
 \hookrightarrow \ldots H^i(\E^{(p^n)}) = 0$$
 for $i > \phi(\E)$ and $n >> 0$.
 
 We can replace $\E$ by $\E^{*}$ and $i$ by $dim X - i$
 in the above sequence of injections. This together with
 Serre duality yields the result.
 \end{proof}
 
When $\E$ is a vector bundle on $\PP^n$, the proposition 
yields the vanishing
$$H^i(\E(j-n-1)) = 0$$
 for $j\ge 0$, $i>\phi(\E)\ge\phi(\E(j))$, 
obtained earlier.

When $k$ is a perfect field, let $W(k)$ be the ring of
Witt vectors over $k$, and $W_2(k)=W(k)/p^2W(k)$.
It is helpful to keep the following example
in mind:  if $k=\Z/p\Z$, then $W(k)$ is the
ring of $p$-adic integers, so that $W_2(k) \cong \Z/(p^2)$.

\begin{thm}\label{thm:van}
Let $k$ be a perfect field of characteristic $p>n$, and
let $X$ be a smooth $n$ dimensional projective $k$-variety
with a reduced normal crossing divisor $D$.
Suppose that  $\E$ is a coherent sheaf on $X$.
If $(X,D)$ can be lifted to a pair over $Spec\, W_2(k)$.
Then
$$H^i(X, \Omega_X^j(\log D)(-D)\otimes \E) = 0$$
for $i+j > n+\phi(\E, D)$.

\end{thm}

\begin{remark}
Note that $\E$ is not required to lift. 
It is possible to obtain a weaker statement  when $p\le n$,
but we won't need it.
\end{remark}

The proof is based on the following lemmas. 

\begin{lemma}\label{lemma:bootstrap1}
 Suppose that $0 \le  D'\le pD$ is a divisor
such that
$$H^i(\Omega_X^j(log D)(-D')\otimes \E^{(p)}) = 0    $$
for all $i+j>N$, then
$$H^i(\Omega_X^j(log D)(-D'_{red})\otimes \E) = 0    $$
for all $i+j>N$.
\end{lemma}

\begin{proof}
Set $D_1 = D'_{red}$ and $B = pD_1-D'$. To avoid confusion,
we will say few words about our conventions. The differentials
on  $\Omega_X^\dt(log D)(B)$ and $\Omega_X^\dt(log D)$ are
inherited from the complex of meromorphic forms.
All other differentials are induced from these using tensor
products and pushforwards.
There is a quasi-isomorphism
$$\Omega_X^\dt(log D)\cong \Omega_X^\dt(log D)(B)$$
where $B = pD_1-D'$ by \cite[3.3]{hara} (see also \cite[4.1]{ms}).
Tensoring both sides with $\E^{(p)}(-D')$ yields
$$\Omega_X^\dt(log D)\otimes \E^{(p)}(-D')\cong
 \Omega_X^\dt(log D)\otimes [\E(-D_1)]^{(p)}.$$
This implies that 
$$F_*(\Omega_X^\dt(log D)(-D')\otimes \E^{(p)})\cong
 F_*(\Omega_X^\dt(log D)\otimes F^*(\E(-D_1))).$$
Lemma~\ref{lemma:projform} shows that the right
side is quasi-isomorphic to 
$$[F_*\Omega_X^\dt(log D)]\otimes \E(-D_1)$$
By  \cite[4.2]{deligne-ill} (and the remarks
of section 1), this is quasi-isomorphic to
$$\left(\bigoplus_j \Omega_X^j(log D)[-j]\right)\otimes \E(-D_1).$$
The spectral sequence
$$H^i(\Omega_X^j(log D)(-D')\otimes \E^{(p)})\Rightarrow
\mathbb{H}^{i+j}(\Omega_X^\dt(log D)(-D')\otimes \E^{(p)})$$
together with the hypothesis
shows that the abutment vanishes for $i+j>N$.
Therefore 
$$\mathbb{H}^{i}(F_*(\Omega_X^\dt(log D)(-D')\otimes \E^{(p)}))
\cong \bigoplus_j H^{i-j}(\Omega_X^j(log D)\otimes  \E(-D_1))$$
vanishes for $i>N$.
 \end{proof}

 \begin{lemma}\label{lemma:bootstrap2}
  Suppose that $0 \le  D'\le p^aD$ is a divisor
 such that
 $$H^i(\Omega_X^j(log D)(-D')\otimes \E^{(p^a)}) = 0    $$
 for all $i+j>N$, then
 $$H^i(\Omega_X^j(log D)(-D'_{red})\otimes \E) = 0    $$
 for all $i+j>N$.
\end{lemma}    

\begin{proof}
We prove this by induction on $a$. 
The case where $D'=0$ is straightforward, so we assume that $D'\not= 0$.
The initial case $a=1$ is the
previous lemma. Suppose that the lemma holds for $a$,
and suppose that $(\E, D')$ satisfies the hypothesis
of the lemma with $a$  replaced by $a+1$.
Since $\{1,\ldots p\}$ forms a set of
representatives of $\Z/p\Z$, we can decompose
$D' = pD_1 + D_2$ such that $D_{2,red} = D'_{red}$, $0< D_2\le pD$ and
$0\le D_1<_{strict} p^aD$. These assumptions guarantee that the divisor
$D'_{red}+D_1$ is less than or equal to $p^aD$ and has the same
support as $D'$.
 Then 
$$\E^{(p^{a+1})}(-D') = \E_1^{(p)}(-D_2)$$
where $\E_1 = \E^{(p^a)}(-D_1)$. By assumption,
$$H^i(\Omega_X^j(log D)(-D')\otimes \E^{(p^{a+1})})
= H^i(\Omega_X^j(log D)(-D_2)\otimes \E_1^{(p)}) = 0$$
for $i+j > N$. Lemma~\ref{lemma:bootstrap1} implies 
$$H^i(\Omega_X^j(log D)(-D'_{red})\otimes \E_1) = 
H^i(\Omega_X^j(log D)(-D'_{red}-D_1)\otimes \E^{(p^{a})})=0$$
for $i+j > N$. Induction yields the desired conclusion.
\end{proof}

\begin{proof} [Proof of theorem \ref{thm:van}]
By definition, $\phi(\E, D) = \phi(\E^{(p^a)}(-D'))$
for some $0\le  D' <_{strict} p^aD$. We can assume $a>>0$ since
we can replace $\E^{(p^a)}(-D')$ by a Frobenius power.
Therefore 
$$H^i(\Omega_X^j(log D)(-p^bD'- D)\otimes \E^{(p^{a+b})}) = 
H^i(\Omega_X^j(log D)(- D)\otimes (\E^{(p^{a})}(-D'))^{(p^b)}) = 0
 $$
for $b>>0$ and all $i>\phi(\E, D) $ and all $j$. 
Since the support of $p^bD'+D$ is $D$ and the coefficients
of this divisor are less than or equal to $p^{a+b}$,
 lemma~\ref{lemma:bootstrap2} implies the theorem.
\end{proof}


\begin{cor}\label{cor:van1}
Suppose that $(X,D,\E)$ is as above and $char\, k = 0$, then
$$H^i(X, \Omega_X^j(log D)(-D)\otimes \E) = 0$$
for $i+j > n + \phi(\E, D)$. In particular,
 $$H^i(X, \Omega_X^j\otimes \E) = 0$$
for $i+j > n+\phi(\E)$
\end{cor}

\begin{proof}
Choose an arithmetic thickening of $(X,D,\E)$. Almost
 all of the closed fibers satisfies the conditions
of the theorem. Therefore the corollary follows by
semicontinuity.
\end{proof}

\begin{cor}\label{cor:van2}
Suppose that $(X,D,\E, k)$ is as in theorem or
corollary~\ref{cor:van1}, and that $L$ is an 
arithmetically nef line bundle  on $X$.
Then
$$H^i(X, \Omega_X^j(log D)(-D)\otimes \E\otimes L) = 0$$
for $i+j > n + \phi(\E, D)$. 
\end{cor}

\begin{proof}
By proposition~\ref{prop:nefline}, $L$ is $F$-semipositive,
hence $\phi(\E\otimes L, D) = \phi(\E, D)$.
So the corollary follows from theorem~\ref{thm:van} in positive characteristic,
or the previous corollary in characterisitic $0$.
\end{proof}

\begin{cor}\label{cor:vandual}
Suppose that $(X,D,\E,k)$ satisfy the conditions of
the theorem or the corollary~\ref{cor:van1}, then
$$Ext^i(\E, \Omega_X^j(log D)) = 0$$
for $i+j < n - \phi(X,D)$. If $\E$ is locally free, then
$$H^i(X,\Omega_X^j(log D)\otimes \E^*) = 0$$
for $i+j < n - \phi(X,D)$.
\end{cor}

\begin{proof} This is a consequence of Serre duality.
\end{proof}

\begin{cor} (Le Potier) 
Suppose that $char \, k = 0$ and $\E_i$ are
ample   locally free sheaves on a smooth variety $X$, then 
  $$H^i(X, \Omega_X^j\otimes \E_1\otimes \ldots \E_m) = 0$$
for $i+j \ge n+rank(\E_1)+\ldots rank(\E_m)$
\end{cor}

\begin{proof} Follows from corollary~\ref{cor:van1} and
  corollary~\ref{cor:keyestimate}.
\end{proof}

The next result is a generalization of the  Kawamata-Viehweg vanishing
theorem \cite{kawamata, viehweg}.
 (To obtain the statement given in the introduction,
set $\Delta = \frac{1}{m}D'$.)

\begin{cor}\label{cor:kvvan1}
With  $(X,D,\E,k)$  as in the corollary~\ref{cor:van1}.
Suppose that there is a positive integer $m$, and
a divisor $0\le D'\le (m-1)D$ such that $S^{m}(\E)(-D')$ is ample.
Then 
$$H^i(X, \Omega_X^j(log D)(-D)\otimes \E) = 0$$
for $i+j\ge n+rank(\E)$. In particular,
$$H^i(X, \omega_X\otimes \E) = 0$$
for $i \ge rank(\E)$.
\end{cor}

\begin{proof}
This follows from theorem~\ref{thm:logkeyestimate} and
corollary~\ref{cor:van1}.
\end{proof}

\begin{cor}\label{cor:genphivan}
Suppose that $char \, k = 0$ and $\E$ is a locally
free sheaf on a projective $k$-variety $Z$ with rational
singularities, then
$$H^i(Z, \omega_Z\otimes \E) = 0$$
for $i> \phi_{gen}(\E)$.
\end{cor}

\begin{proof}
By the previous corollary, $H^i(Y,\omega_Y\otimes f^*\E)$
vanishes for  $i> \phi_{gen}(\E)$, for some resolution
of singularities $f:Y\to Z$. We have $\R f_*\omega_Y = \omega_Z$
 because $Z$ has rational singularities, and so the
corollary follows.
\end{proof}

\begin{cor}\label{cor:kvvanishing} 
Suppose that $char \, k = 0$ and $\E$ is the
pull back of an ample vector bundle
under a surjective morphism $f:X\to Y$ with $X$
smooth, then
$$H^i(X, \omega_X\otimes \E) = 0$$  
for $i\ge rank(\E)+d$ where $d$ is  the dimension
of the generic fiber of $f$.
\end{cor}   

\begin{proof}
Follows from 
corollary~\ref{cor:genphi} and corollary~\ref{cor:genphivan}.
\end{proof}

\section{Some applications}

In this section, we work over $\C$.
We start with a  refinement of the Lefschetz
hyperplane  theorem (which corresponds to  case B with 
$\E= O_X$ and $r=m=1$).

\begin{prop}
 Let $D\subset X$ be a smooth divisor on a smooth $n$ dimensional 
projective variety.
Suppose that $m>0$.
\begin{enumerate}
\item[A.] If   $S^m(\E)(rD)$ is ample for some $-m< r\le 0$, then 
$$H^i(X,\Omega_X^j\otimes \E)\to H^i(D,\Omega_D^j\otimes \E)$$
is bijective if $i+j\ge n + rank(\E)$ and surjective if
$i+j= n +rank(\E)-1$.

\item[B.] If $S^m(\E)(rD)$ is ample for some $0< r \le m$, then 
$$H^i(X,\Omega_X^j\otimes \E^*)\to H^i(D,\Omega_D^j\otimes \E^*)$$
is bijective if $i+j\le n - rank(\E)-1$ and injective if
$i+j= n -rank(\E)$.
\end{enumerate}

\end{prop}

\begin{proof}
We have an exact sequence
$$0\to \Omega_X^j(log D)(-D)\to \Omega_X^j\to \Omega_D^j\to 0.$$
Tensoring this with $\E$ and applying  corollary  \ref{cor:kvvan1}
proves A. For B, tensor this with $\E^*$ and observe that
by Serre duality and corollary \ref{cor:kvvan1}
$$H^i(X,\Omega_X^j(log D)(-D)\otimes \E^*)\cong
H^{n-i}(X,\Omega_X^{n-j}(log D)(-D)\otimes \E(D))^*=0 $$
when $i+j \le n-rank(\E)$.
\end{proof}

 Given an algebraic variety
$Y$ and algebraic coherent sheaf $\F$ over $Y$, we denote
the corresponding analytic objects by $Y^{an}$ and
$\F^{an}$. Given a closed subvariety $Z\subset Y$, let
$$codim(Z,\,Y) = 
\left\{
\begin{array}{l}
 \mbox{ the codimension of $Z$ if $Z\not=\emptyset$} \\
 $dim\, Y + 1$\mbox{ otherwise}
\end{array}
\right.
$$

Our goal  is to prove a vanishing theorem
for ample vector bundles over quasiprojective  varieties.
This generalizes some results for line bundles due
to  Bauer, Kosarew \cite{bk2} and the author \cite{arapura1}. 

\begin{thm}\label{thm:noncompact}
Suppose that $U$ is a smooth quasiprojective variety
with a possibly singular projective compactification
$Y$.  Let $\E$ be the restriction to $U$ of an ample
vector bundle on some compactification of $U$ (possibly other than
$Y$). Then  
$$H^i(U^{an},\, (\Omega_U^j\otimes \E^*)^{an})= 
H^i(U,\, \Omega_U^j\otimes \E^*) = 0$$
for $i+j\le  codim(Y-U,\, Y)- rank(\E)-1$.
\end{thm}

As a first step, we need a generalization of  Steenbrink's vanishing
theorem \cite{steenbrink}.

\begin{prop}\label{prop:steen}
Suppose that $f:X\to Y$ is a desingularization of an $n$-dimensional
 projective variety such that $X$ possesses a reduced normal crossing
divisor $D$ containing the exceptional divisor.
 If $\E$ is a nef (e. g. globally generated) vector
bundle, then 
$$R^if_*[\Omega_X^j(log D)(-D)\otimes \E] = 0$$
for $i+j \ge n + rank(\E)$.
\end{prop}

\begin{proof}
As in the proof of corollary~\ref{cor:biratphi}, we can find
a relatively ample divisor $-\sum a_iD_i$ supported on $D$ 
with $a_i\ge 0$. Let $L$ be a large power of an ample line
bundle on $Y$. Then
$f^*L^{\otimes N}(-\sum a_iD_i)$ is ample for  all $N>0$.
Therefore  $S^m(\E\otimes f^*L^{\otimes N})(-\sum a_iD_i)$ is ample for all
$m> 0$. Then corollary~\ref{cor:kvvan1} shows that 
$$H^i(\Omega_X^j(log D)(-D)\otimes \E\otimes f^*L^{\otimes N}) = 0$$
for $i+j \ge n + rank(\E)$.
For $N>>0$, the Leray spectral sequence and Serre vanishing yields
$$H^0(R^if_*[\Omega_X^j(log D)(-D)\otimes \E]\otimes L^{\otimes N})
= 0$$
for $i+j \ge n + rank(\E)$.
We can assume that these sheaves are all globally generated by
increasing $N$ if necessary. Thus they must vanish.
\end{proof}

\begin{cor}\label{cor:localvan}
Let $Y$ be a  projective variety,
$Z\subset Y$ a closed subvariety containing the
singular locus, and  $f:X\to Y$ a resolution of singularities
which is an isomorphism over $Y-Z$ 
such that $D= f^{-1}(Z)_{red}$ is a divisor with normal crossings.
Then for any nef vector bundle $\E$ on $X$,
$$H_D^i(X,\Omega_X^j(log D)\otimes \E^*) = 0$$
for $i+j\le codim(Y,Z)  -rank(\E)$.
\end{cor}

\begin{proof} This is a generalization of 
\cite[thm 1]{arapura1}. A proof of this corollary
can be obtained by simply replacing Steenbrink's theorem
with proposition~\ref{prop:steen} in the proof given there.
\end{proof}

\begin{proof}[Proof of theorem~\ref{thm:noncompact}]
The groups $H^i(U, \F)$ and  $ H^i(U^{an}, \F^{an})$
are isomorphic for any coherent sheaf $\F$ on $Y$
and $i< codim(Y-U,Y)-1$ by \cite[IV 2.1]{hartshorne-ampsub}.
So it suffices to prove the vanishing in the algebraic category.
Let $Z = Y-U$ and let $f:X\to Y$ be a desingularization
satisfying the assumptions of corollary~\ref{cor:localvan}.
We can assume that $X$ also dominates the compactification
where $\E$ extends to an ample bundle. We use the same
symbol for this extension, and its pullback to $X$.
Corollary~\ref{cor:biratphi}
and theorem~\ref{thm:keyestimate} implies that
$\phi(\E, D) \ < rank(\E)$.
Corollary~\ref{cor:localvan} shows that 
$$H^{i+1}_D(X,\Omega_X^j(log D)\otimes \E^*)=0.$$
Then the exact sequence for local cohomology yields
a surjection
$$H^i(X,\Omega_X^j(log D)\otimes \E^*)\to
H^i(U,\Omega_X^j(log D)\otimes \E^*) =
H^i(U,\Omega_U^j\otimes \E^*).$$
The cohomology group on the left
vanishes as a consequence of  corollary~\ref{cor:vandual}.
\end{proof}


\newpage

\begin{appendix}
\section{Arithmetically nef bundles}
\centerline{{Dennis S. Keeler\footnote{ 
        Partially supported by an NSF Postdoctoral Research Fellowship. \\
        Department of Mathematics, MIT, Cambridge, MA 02139, \texttt{dskeeler@mit.edu}}}}
        
Let $Y$ be a noetherian scheme, let
$f\colon X \to Y$ be a proper morphism, and let $\E$ be a vector bundle on $X$.
For each $y \in Y$, let $\E_y$ be the restriction of $\E$ to the fiber $X_y$.
Recall that $\E$ is \emph{$f$-nef} if $\E_y$ is nef for every closed
$y \in Y$ (see, for instance, \cite[Definition~2.9]{keeler}).
If $Y$ is affine, then the property of $\E$ being $f$-nef does not
depend on $f$, so we may simply say that $\E$ is nef \cite[Proposition~2.15]{keeler}.
        
\begin{defn}\label{def:arithmetically-nef}
Let $X$ be a proper scheme over a field $k$, and let $\E$ be a vector bundle on $X$.
If $char\, k = 0$, then $\E$ is \emph{arithmetically nef} if there exists a thickening
$(\tilde{X} \to Spec\, A, \tilde{\E})$ such that $\tilde{\E}$ is nef.
For convenience, if $char\, k = p > 0$, we say that $\E$ is arithmetically nef
if $\E$ is nef.
\end{defn}

Note that if $\E$ is arithmetically nef, then $\tilde{\E}$ will be nef
on every fiber of a certain thickening \cite[Lemma~2.18]{keeler}. In particular,
if $\E$ is arithmetically nef, then $\E$ is nef.

Like nefness, the property of being arithmetically nef behaves well under
pullbacks.

\begin{lemma}\label{lem:arith-nef-surjective}
Let $\E$ be a vector bundle on a proper scheme $X$ over a field $k$, and
let $f\colon X' \to X$ be a proper morphism. 
\begin{enumerate}
\item If $\E$ is arithmetically nef, then $f^*\E$ is arithmetically nef, and
\item if $f$ is surjective and $f^*\E$ is arithmetically nef, then
$\E$ is arithmetically nef.
\end{enumerate}
\end{lemma}
\begin{proof}
By replacing $X$ with $\PP(\E)$ and $X'$ with $\PP(f^*\E)$, we may assume that
$\E$ equals a line bundle $L$. 

If $char\, k = 0$, then choose a thickening
$(\tilde{f}\colon \tilde{X'} \to \tilde{X}, \tilde{L})$ such
that $\tilde{f}$ is proper and also surjective if $f$ is surjective
\cite[$\mathrm{IV}_3$, 8.10.5]{ega-app}. If $L$ is arithmetically nef, then upon
further shrinking the thickening 
we may assume $\tilde{L}$ is nef,
and so $\tilde{f}^* \tilde{L}$ is nef \cite[Lemma~2.17]{keeler}.
Hence $f^*L$ is arithmetically nef. Now if $f^*L$ is arithmetically nef,
then $L$ is arithmetically nef by a similar argument \cite[loc. cit.]{keeler}.
If $char\, k = p > 0$, then these are just statements about nef line bundles
\cite[loc. cit.]{keeler}.
\end{proof}

Being arithmetically nef also behaves well under tensor product.

\begin{lemma}\label{lem:arith-nef-tensor}
Let $L$ and $M$ be line bundles on a proper scheme $X$. Then
\begin{enumerate}
\item $L$ is arithmetically nef, if and only if $L^n$ is arithmetically nef for
all $n > 0$, if and only if $L^n$ is arithmetically nef for some $n > 0$, and
\item If $L$ and $M$ are arithmetically nef, then $L \otimes M$ is arithmetically nef.
\end{enumerate}
\end{lemma}
\begin{proof} These statements follow immediately from the definition of nef.
\end{proof}

It is natural to conjecture that if $L$ is nef, then $L$ is arithmetically nef, but
we have been unable to prove this. However, we do have a non-trivial collection of examples.

\begin{prop}\label{prop:arith-nef-examples}
Let $L$ be a line bundle on a proper scheme $X$ over a field $k$. If $L$ is 
 semi-ample (i.e., there exists $n > 0$ such that $L^n$ is generated by global sections)
 or $L$ is numerically trivial,
then $L$ is arithmetically nef.
\end{prop}
\begin{proof}
The case $char\, k = p>0$ is trivial, so assume $char\, k = 0$.
First suppose that $L$ is semi-ample. We may replace $L$ with $L^n$ \eqref{lem:arith-nef-tensor}
 and assume
that $L$ is generated by global sections. Then $L$ defines a $k$-morphism $f\colon X \to \PP^m$
for some $m$, and $L = f^*O_{\PP^m}(1)$ \cite[II 7.1]{hartshorne2-app}. Now $f$ is
proper \cite[II 4.8e]{hartshorne2-app}, so we may replace $L$ with $O_{\PP^m}(1)$ by
lemma~\ref{lem:arith-nef-surjective}. Now since $O(1)$ is ample, there exists
a thickening such that $\widetilde{O_{\PP^m}(1)}$ is ample \cite[$\mathrm{III}_1$, 4.7.1]{ega-app}.
Hence $O_{\PP^m}(1)$ is arithmetically nef.

Now suppose that $L$ is numerically trivial. Any pullback of $L$ is also
numerically trivial \cite[Lemma~2.17]{keeler}. Thus by lemma~\ref{lem:arith-nef-surjective}
we may replace $X$ with a Chow cover and thus assume that $X$ is projective.
We may also replace $X$ with the disjoint union
$\coprod_i X_i$, where the $X_i$ are the reduced, irreducible
components of $X$. Thus we may assume that $X$ is integral.  There exists
a projective, surjective morphism $X' \to X$ such that $X'$ is geometrically
integral \cite[Lemma~3.3]{keeler}, and thus we may assume that $X$ is geometrically integral.

Let $H$ be an ample divisor on $X$. We may choose a thickening 
$(\tilde{\pi}\colon\tilde{X} \to Spec\, A, \tilde{H}, \tilde{L})$ such that $\tilde{H}$
is ample \cite[$\mathrm{III}_1$, 4.7.1]{ega-app}
and $\tilde{\pi}$ is flat \cite[$\mathrm{IV}_3$, 8.9.4]{ega-app}. Further,
we may assume that all fibers of $\tilde{\pi}$ are geometrically integral
\cite[$\mathrm{IV}_3$, 12.2.4]{ega-app}.

For any line bundles $N, M$ on $\tilde{X}$ and $s \in Spec\, A$, let $N_s, M_s$ be the
restriction of $N, M$ to the fiber $\tilde{X}_s$.
Then since $\tilde{\pi}$ is flat, the intersection numbers
$(N_s^r.M_s^{\dim X - r})$ are independent of $s$ \cite[Remark~3.5]{keeler}. Now
$\tilde{L}_s$ is numerically trivial if and only if
\[
(\tilde{L}_s.\tilde{H}_s^{\dim X -1}) = (\tilde{L}_s^2.\tilde{H}_s^{\dim X - 2}) =  0
\]
by \cite[p.~305, Corollary~3]{kleiman}. But since these intersection numbers
are $0$ at the generic point, they are $0$ at each $s \in Spec\, A$. 
Thus $L$ is arithmetically nef.
\end{proof}

\begin{cor}\label{cor:arith-nef-on-curves}
Let $X$ be a projective scheme with $\dim X \leq 1$. If $L$ is a nef line bundle,
then $L$ is arithmetically nef.
\end{cor}
\begin{proof} Using lemma~\ref{lem:arith-nef-surjective}, we may assume that
$X$ is integral. Then $L$ is either numerically trivial
or ample, and hence arithmetically nef by proposition~\ref{prop:arith-nef-examples}.
\end{proof}

We now consider arithmetically nef line bundles on a surface $X$, that is, an
integral scheme of dimension $2$.

\begin{cor}\label{cor:arith-nef-on-surface}
Let $X$ be a projective surface, and let $L$ be a nef line bundle
such that $L^n$ is effective for some $n > 0$ (for example, if $L$ is big). Then $L$ is
arithmetically nef.
\end{cor}
\begin{proof}
We assume that the characteristic of the ground field is $0$. Let $H$ be an ample line
bundle.
By lemma~\ref{lem:arith-nef-tensor} we may replace $L$ by $L^n$ and hence assume
that $L \cong O(D)$ for an effective divisor $D$.
Since $L\vert_D$ is arithmetically nef by corollary~\ref{cor:arith-nef-on-curves},
we may choose an arithmetic thickening 
$(\tilde{X}, \tilde{H},\tilde{L} \cong \tilde{O}(\tilde{D}))$
such that $\tilde{L}\vert_{\tilde{D}}$ is nef
and $\tilde{H}$ is ample \cite[$\mathrm{III}_1$, 4.7.1]{ega-app}.

Consider the short exact sequences
\[
H^i(\tilde{X}, \tilde{H}^n \otimes \tilde{L}^m)
\to H^i(\tilde{X}, \tilde{H}^n \otimes \tilde{L}^{m+1})
\to  H^i(\tilde{D}, \tilde{H}\vert_{\tilde{D}}^n \otimes \tilde{L}\vert_{\tilde{D}}^{m+1})
\]
with $i > 0, m\geq 0, n > 0$. We may fix $n$ sufficiently large so that the leftmost
group vanishes for $m = 0$ and the rightmost vanishes for all $m \geq 0$ 
\cite[Theorem~1.5]{keeler}. But then by induction on $m$, the middle vanishes for all $m \geq 0$.
Since any coherent sheaf $\F$ on $\tilde{X}$ is a quotient of a finite direct sum
of $\tilde{H}^{-\ell}$, we have that for any $\F$, there exists $n$ such
that $H^i(\tilde{X}, \F \otimes \tilde{H}^{n} \otimes \tilde{L}^m) = 0$ for $i > 0, m \geq 0$,
and so $\tilde{L}$ is nef 
\cite[Proposition~5.18]{keeler}. Thus $L$ is arithmetically nef.
\end{proof}

\section{Arithmetically nef line bundles are $F$-semipositive}

In this section, we characterize $F$-semipositive line bundles. Lemma~\ref{lemma:ampleline}
states that a line bundle $L$ is $F$-ample if and only if it is ample.
We later see in lemma~\ref{lemma:Fsemiposisnef} 
that any $F$-semipositive
vector bundle is arithmetically nef. Given these two facts, it is natural to conjecture
that a line bundle $L$ is $F$-semipositive if and only if it is arithmetically
nef, and indeed this
is the case.

\begin{prop}\label{prop:nefline}
Let $X$ be a projective variety over a field $k$ and let $L$ be a line bundle
on $X$. Then $L$ is $F$-semipositive if and only if $L$ is arithmetically nef.
\end{prop}
\begin{proof}
Given lemma~\ref{lemma:Fsemiposisnef}, we need only show that an arithmetically nef 
line bundle is $F$-semipositive.
Let $O_X(1)$ be a very ample line bundle.
If $char\, k =p >0$, there exists $m$ such that
$H^i(X, O_X(m-i) \otimes L^n) = 0$ for $i>0, n \geq 0$ \cite[Theorem~1.5]{keeler}.
Thus $areg(L) \leq m$, so $L$ is $F$-semipositive.

If $char\, k = 0$, then choose a thickening $(\tilde X, {\tilde O_X(1)}, \tilde L)$
over a finitely generated $\Z$-algebra $A \subset k$ such that ${\tilde O_X(1)}$ is very ample
\cite[$\mathrm{III}_1$, 4.7.1]{ega-app}
and $\tilde L$ is nef. Again there exists $m$ such that
$H^i({\tilde X}, {\tilde O_X(m-i)} \otimes {\tilde L}^n) = 0$ for $i>0, n \geq 0$ 
\cite[loc. cit.]{keeler}.
So by semicontinuity, for each closed point $q$ of $Spec\, A$, we
have $areg(L_q) \leq m$. Thus $areg(L) \leq m$ and $L$ is again $F$-semipositive.
\end{proof}

\section{Base change}
Most of the seminal works on ample vector bundles, such as \cite{hartshorne-app,barton,gieseker-app},
assumed that the base field $k$
was algebraically closed and some of their proofs use this assumption.
However, since this paper's ``reduction to characteristic $p$'' methods require the
non-algebraically closed case, we now prove a few standard lemmas which will allow
the application of
``algebraically closed results'' to the general case.

Since we do not assume that our varieties are geometrically integral, we must also allow
$X$ to be a projective scheme. For a coherent sheaf $\F$ we keep the same definitions of
$\phi(\F)$, $F$-ample, and $F$-semipositive as given in sections \ref{sec:frobenius-amplitude} 
and \ref{sec:asymptotic-regularity}, 
just with $X$ as a scheme, projective over a field $k$. 

\begin{lemma}\label{lem:reduced}
Let $X$ be a projective scheme over a field $k$ of characteristic $p > 0$,
let $\E$ be a vector bundle on $X$, and let $\E_\red$
be the restriction of $\E$ to the reduced scheme $X_\red$. 
Then $\phi(\E) =  \phi(\E_\red)$.
\end{lemma}
\begin{proof}
This follows from a standard argument as 
in \cite[ex. III 3.1, 5.7]{hartshorne2-app},
using the commutative diagram \eqref{eq:frobenius-diagram} 
with $Y = X_\red$. To see this, let $f\colon X_\red \to X$
be the natural immersion, let $\calN$ be the nilradical of $O_X$, and let
$\F$ be a coherent sheaf on $X$. Then $\calN^{r_0} = 0$ for $r_0$ sufficiently large
and we have short exact sequences
\[
0 \to \calN^{r+1}\F \to \calN^r\F \to \calN^r\F/\calN^{r+1}\F \to 0
\]
with $0 \leq r < r_0$.
Since $f_*O_{X_\red} = O_X/\calN$, there exist coherent $\G_r$ on $X_\red$
such that $f_* \G = \calN^r\F/\calN^{r+1}\F$.
Thus, if $i > \phi(\E_\red)$, we have
\[
H^i(X_\red, \G_r \otimes f^*\E^{(p^n)}) = H^i(X, \calN^r\F/\calN^{r+1}\F \otimes \E^{(p^n)})=0
\]
for $0 \leq r < r_0$ and $n \gg 0$.
Descending induction on $r$ then gives $H^i(X, \F \otimes \E^{(p^n)})=0$ for $n \gg 0$.
\end{proof}

\begin{lemma}\label{lem:irreduciblecomponents}
Let $X$ be a reduced projective scheme over a field $k$ of characteristic $p > 0$,
let $X_i$ be the irreducible components of $X$, $i=1,\dots,r$,
 let $\E$ be a vector bundle on $X$, and let $\E_i$ be the restriction
 of $\E$ to $X_i$. Then $\phi(\E) = \max_i \phi(\E_i)$.
\end{lemma}
\begin{proof}
This follows from a standard argument as 
in \cite[ex. III 3.2, 5.7]{hartshorne2-app},
using the commutative diagrams \eqref{eq:frobenius-diagram} 
for each $i = 1, \dots, r$ with $Y = X_i$.
Let $\F$ be a coherent sheaf on $X$ and let $\calI_j$ be the sheaf of ideals of $X_j$.
Then $\calI_j\F$ has support on the $X_i$ with $i \neq j$ and $\F/\calI_j\F$ has
a natural $O_{X_j}$-module structure. So 
by induction on the number of irreducible components $r$, 
for $i > \max_i \phi(\E_i) \geq \phi(\E_j)$, we have
short exact sequences
\[
0 = H^i(\calI_j \F \otimes \E^{(p^n)}) \to H^i(\F \otimes \E^{(p^n)})
\to H^i(\F/\calI_j\F \otimes \E^{(p^n)}) = 0
\]
for $n \gg 0$. Thus $\phi(\E) \leq \max_i \phi(\E_i)$, and the reverse inequality is trivial.
\end{proof}

\begin{lemma}\label{lem:finite-surjective}
Let $f\colon Y \to X$ be a finite, surjective morphism (not necessarily a $k$-morphism)
of projective $k$-schemes, and
let $\E$ be a vector bundle on $X$. Then $\phi(\E) = \phi(f^*\E)$.
\end{lemma}
\begin{proof}
When $char\, k = p > 0$, this again follows from a standard argument as 
in \cite[ex. III 4.2, 5.7]{hartshorne2-app},
using the commutative diagram 
\eqref{eq:frobenius-diagram}. 
By lemmas~\ref{lem:reduced} and \ref{lem:irreduciblecomponents}, we may assume
that $X$ and $Y$ are integral schemes. Then one may follow the argument
outlined in \cite[ex. III 4.2]{hartshorne2-app}.

When $char\, k = 0$, we can choose a thickening 
$({\tilde f}\colon {\tilde Y} \to {\tilde X}, {\tilde \E})$ 
and assume that ${\tilde f}$ is finite and surjective 
\cite[$\mathrm{IV}_3$, 8.10.5]{ega-app}. The claim
then follows from the positive characteristic case.
\end{proof}

\begin{cor}\label{cor:reduced-irreducible-char-0}
Let $X$ be a projective scheme over a field of arbitrary characteristic. Then
the conclusions of lemmas~\ref{lem:reduced} and \ref{lem:irreduciblecomponents}
remain true.
\end{cor}
\begin{proof} 
Let $\coprod_i X_i$ be the disjoint union of the irreducible components of $X$ if $X$ is reduced.
Both claims follow from lemma~\ref{lem:finite-surjective} because
the maps
$X_\red \to X$ and $\coprod_i X_i \to X$ are finite and surjective.
\end{proof}

We now show that all of our concepts of ampleness behave well under base change.

\begin{lemma}\label{lem:basechange}
Let $k \subseteq k'$ be fields, let $X$ be a projective scheme over $k$, and let $\E$
be a vector bundle on $X$. Then $\phi(\E) = \phi(\E \otimes_k k')$. Also,
 each of the following properties hold for $\E$ if and only if they hold for $\E \otimes_k k'$
 on $X \times_k k'$:
 \begin{enumerate}
 \item\label{item:F-ample} $F$-ample,
 \item\label{item:F-semipositive} $F$-semipositive,
 \item\label{item:p-ample} $p$-ample, if $char\, k = p > 0$,
 \item\label{item:ample} ample, 
 \item\label{item:nef} nef,
 \item\label{item:arith-nef} arithmetically nef.
 \end{enumerate}
\end{lemma}
\begin{proof}
Let $f\colon X \times_k k' \to X$ be the base change, and
let $O_X(1)$ be a very ample line bundle for $X$. Then $f^*O_X(1)$ is very ample on $X \times k'$
\cite[II, 4.4.10]{ega-app}. If $char\, k = p > 0$, then we have a commutative diagram
\eqref{eq:frobenius-diagram} with $Y = X \times k'$. Note that this
diagram is \emph{not} cartesian because $F_{Y} \neq F_X \times {\mathrm id}_{k'}$.
However, the commutivity of \eqref{eq:frobenius-diagram} gives
\begin{multline*}
H^i(Y, ( f^* \E)^{(p^n)} \otimes f^* O_X(b) ) =
H^i(Y,  F^{n*}_{Y} f^* \E \otimes f^* O_X(b) ) \\ =
H^i(Y, f^* F^{n*}_X \E \otimes f^* O_X(b) ) = 
H^i(Y, f^* (\E^{(p^n)}) \otimes f^* O_X(b) )
\end{multline*}
for $i \geq 0, n\geq 0, b \in \Z$.
Now since $k \to k'$ is a flat morphism \cite[III 9.3]{hartshorne2-app},
\[
 H^i(Y, (f^*\E)^{(p^n)} \otimes f^*O_X(b)) = H^i(X, \E^{(p^n)} \otimes O_X(b)) \otimes_k k'.
\]
Then $\phi(\E) = \phi(f^*\E)$ follows by corollary~\ref{cor:FamplViaLineb}
since $k \to k'$ is faithfully flat. If $char\, k = 0$, it is clear by the definition
of $\phi$ that $\phi(\E) = \phi(f^*\E)$ because an arithmetic thickening of $X$
is an arithmetic thickening of $X \times k'$. Now \eqref{item:F-ample} is immediate
and a similar proof gives \eqref{item:F-semipositive}.

For \eqref{item:p-ample}, consider the exact sequence
\[
H^0(X, \E^{(p^n)} \otimes O_X(b)) \otimes_k O_X \to \E^{(p^n)} \otimes O_X(b) \to \calC_n \to 0
\]
for fixed $b \in \Z$.
Since $k \to k'$ is faithfully flat, we obtain a similar map of global sections
of $(f^*\E)^{(p^n)} \otimes f^*O_X(b)$ by tensoring the exact sequence with $- \otimes_k k'$
(or equivalently, pulling back by $f$).
Then the cokernel $\calC_n = 0$ if and only if $\calC_n \otimes_k k' = 0$, so $\E$ is $p$-ample
if and only if $f^*\E$ is $p$-ample.

Finally, \eqref{item:ample}--\eqref{item:arith-nef} can be reduced to the case of a line bundle
by working on $\PP(\E)$ and $\PP(f^*\E)$. The cases \eqref{item:ample} and
\eqref{item:nef} are given by 
\cite[$\mathrm{IV}_2$, 2.7.2]{ega-app} and \cite[Lemma~2.18]{keeler}. 
 If $char\,k = 0$,
then \eqref{item:arith-nef} is clear because an arithmetic thickening of $X$
is an arithmetic thickening of $X \times k'$.
\end{proof}

\section{Dimensional bound on $F$-amplitude of ample bundles}
Let $\E$ be an ample vector bundle. We have seen that if $char\, k = 0$, then
$\phi(\E) < rank(\E)$ (theorem~\ref{thm:keyestimate}). We will now derive
another bound on $\phi(\E)$, which is independent of the characteristic of $k$.
First, we need some lemmas.

\begin{lemma}\label{lem:char-p-very-ample}
Let $X$ be a projective scheme over a field $k$ of characteristic $p > 0$,
let $O_X(1)$ be a very ample invertible sheaf, and let $\E$ be a vector bundle. 
Then for any $b \in \Z$, there exists $n_0$ such that
\[
H^i(O_X(b+m) \otimes \E^{(p^n)}) = 0
\]
for all $i > \phi(\E), n \geq n_0, m \geq 0$.
\end{lemma}
\begin{proof}
We induct on $dim\, X$; the case of $dim\, X = 0$ is trivial. Since $O_X(1)$ is 
very ample, we may choose an effective Cartier divisor $H$ with $O_X(H) \cong O_X(1)$
and $dim\, H = dim\, X - 1$. By theorem~\ref{thm:first} \eqref{thm:first4} 
generalized to the case of schemes,
$\phi(\E_H) \leq \phi(\E)$. So the claim follows from the short exact sequences
\[
H^i(O_X(b+m) \otimes \E^{(p^n)}) \to H^i(O_X(b+m+1) \otimes \E^{(p^n)}) 
\to H^i(O_H(b+m+1) \otimes \E_H^{(p^n)}),
\]
where the case $m = 0$ follows from the fact that $\phi(\E) < i$.
\end{proof}

\begin{lemma}\label{lem:very-ample-bounds}
Let $X$ be a projective scheme, let $H$ be a very ample Cartier divisor,
and let $\E$ be a vector bundle. Then
\[
\phi(\E_H) \leq \phi(\E) \leq \phi(\E_H) + 1.
\]
\end{lemma}
\begin{proof}
The first inequality is just theorem~\ref{thm:first} \eqref{thm:first4},
generalized to schemes.

If $char\, k = p > 0$, then applying lemma~\ref{lem:char-p-very-ample} to
$H$, for any $b \in \Z$ there exists $n_0$ such that there exist exact sequences
\begin{multline*}
0 = H^i(O_H(b+m+1) \otimes \E_H^{(p^n)}) 
\to H^{i+1}(O_X(b+m) \otimes \E^{(p^n)}) 
\\ \to H^{i+1}(O_X(b+m+1) \otimes \E^{(p^n)})
\to H^{i+1}(O_H(b+m+1) \otimes \E_H^{(p^n)}) = 0
\end{multline*}
for $i > \phi(\E_H), n \geq n_0, m \geq 0$. By Serre Vanishing, 
$H^{i+1}(O_X(b+m+1) \otimes \E^{(p^n)}) = 0$ for $m \gg 0$. So by descending induction on $m$
and corollary~\ref{cor:FamplViaLineb},
$\phi(\E) \leq \phi(\E_H) + 1$.
The case of $char\, k = 0$ is then immediate.
\end{proof}

It is now an easy matter to obtain our bound on $\phi(\E)$ for ample $\E$.
This generalizes lemma~\ref{lemma:phipample} and proposition~\ref{prop:curveFample-ample}. 

\begin{prop}\label{prop:ample-dimension-bound}
Let $X$ be a projective scheme over a field $k$ with
 $dim\, X > 0$, and let $\E$ be an ample vector bundle.
Then $\phi(\E) < dim\, X$.
\end{prop}
\begin{proof}
We may assume that $k$ is algebraically closed \eqref{lem:basechange} and that
$X$ is reduced, irreducible \eqref{cor:reduced-irreducible-char-0}, and normal
\eqref{lem:finite-surjective}.
If $dim\, X = 1$, then the claim is proposition~\ref{prop:curveFample-ample}.
If $dim\, X > 1$, then induction on $dim\, X$ and lemma~\ref{lem:very-ample-bounds}
yields the result.
\end{proof}


\begin{thebibliography}{ABC}

\bibitem[A1]{arapura1} D. Arapura, {\em Local cohomology of sheaves of
    differential forms and Hodge theory}, J.  Reine
  Angew. Math. 409, (1990)

\bibitem[A2]{arapura} D. Arapura, {\em A class of
sheaves satisfying Kodaira's vanishing theorem},
Math. Ann. 318 (2000)

\bibitem[BK1]{bk} I. Bauer, S. Kosarew, {\em On the Hodge spectral
sequence for some classes of non-complete algebraic varieties},
Math. Ann. 284 (1989)

\bibitem[BK2]{bk2} I. Bauer, S. Kosarew, {\em Some aspects of
   Hodge theory on noncomplete algebraic manifolds} Prospects in complex
   geometry (Katata and Kyoto, 1989), 281--316, Lect. Notes in Math.,
   1468, Springer-Verlag (1991).

\bibitem[CL]{carter-lustig} R. Carter, G. Lusztig,
{On modular representations of the general linear
group and symmetric groups}, Math. Z. 136 (1974), 193- 242

\bibitem[C]{Cataldo} M. de Cataldo, {\em Vanishing
via lifting to second Witt vectors and a proof of
an isotriviality result}, J. Algebra 219 (1999), 255-265

\bibitem[DI]{deligne-ill} P. Deligne, L. Illusie,
{\em Relevetments modulo $p^2$  et decomposition du
complexe de de Rham}, Inv. Math. 89 (1987), 247-280

\bibitem[EGA]{ega} A. Grothendieck, J. Dieudonn\'e,
{\'El\'ements de g\'eom\'etrie alg\'ebrique}
Publ. IHES (1960-1967)
 
\bibitem[EV]{esnault-v} H. Esnault, E. Viehweg, {\em
Lectures on vanishing theorems}, Birkh\"auser (1993)

\bibitem[F]{fulton} W. Fulton, {\em Young Tableaux}
Cambridge U. Press (1997)

\bibitem[FL]{fulton-lang} W. Fulton, S. Lang, {\em Riemann-Roch
Algebra}, Springer-Verlag (1985)

\bibitem[G]{gieseker} D. Gieseker, {\em P-ample bundles and their 
Chern classes}, Nagoya Math. J 43 (1971), 91-116

\bibitem[Ha]{hara} N. Hara, {\em A characterization of 
rational singularities in terms of injectivity of Frobenius maps}
Amer. J. Math. 120 (1998), 981-996

\bibitem[H1]{hartshorne} R. Hartshorne, {\em Ample vector bundles},
  Publ. IHES 29 (1966), 63-94

\bibitem[H2]{hartshorne-ampsub}  R. Hartshorne, {\em Ample
    subvarieties of algebraic varieties}, Lect. Notes in Math. 156,
Springer-Verlag (1970)

\bibitem[H3]{hartshorne2} R. Hartshorne, {\em Algebraic geometry},
Springer-Verlag (1977)

\bibitem[Ho]{horrocks} G. Horrocks, {\em Vector bundles on the
punctured spectrum of a local ring } Proc. Lond. Math. Soc. 14
(1964), 689-713

\bibitem[I]{illusie} L. Illusie, {\em R\'eduction semi-stable
et d\'ecomposition de complexes de de Rham \`a coefficients}
Duke Math. J. 60 (1990)

\bibitem[I2]{illusie2} L. Illusie, {\em Frobenius et
d\'eg\'en\'erescence de Hodge}
Introduction \`a le Th\'eorie de Hodge, Soc. Math. Frances (1996)

\bibitem[J]{jantzen} J. C. Jantzen, {\em Representations
of algebraic groups},  Academic Press (1987)

\bibitem[Ka]{kawamata} Y. Kawamata, {\em A generalization of 
Kodaira-Ramanujam's vanishing theorem}, Math. Ann. 261 (1982) 43-46


\bibitem[L]{lepotier} J. Le Potier, {\em Annulation de la
cohomolgie \'a valeurs dans un fibr\'e vectorial
holomorphe positif de rang quelconque}
Math. Ann. 218 (1975)



\bibitem[MR]{mehta} V. Mehta, A. Ramanathan, {\em Frobenius
splitting and cohomology vanishing for Schubert varieties}
Ann. Math 122 (1985)

\bibitem[MS]{ms} V. Mehta, V. Srinivas, {\em
A characterization of rational singularities}, Asian J. Math. 1
(1997), 249-271

\bibitem[Mi]{mi} L. Migliorini, {\em Some observations on
cohomologically $p$-ample bundles.}  
Ann. Mat. Pura Appl.   164  (1993), 89--102

\bibitem[M]{mumford} D. Mumford, {\em Lectures on curves on an algebraic
surface.} Princeton Univ. Press (1966)

\bibitem[SS]{sh-sommese} B. Shiffman, A. Sommese, {\em Vanishing
theorems on complex manifolds}, Birk\"auser (1985)

\bibitem[S]{smith} K. Smith, {\em Vanishing, singularities and
effective bounds via prime characteristic local algebra},
Algebraic Geometry, Santa Cruz 1995, AMS. 

\bibitem[St]{steenbrink} J. Steenbrink, {\em Vanishing theorems
for singular spaces}, Ast\'erisque 130 (1984), 330-341

\bibitem[V]{viehweg} E. Viehweg,{\em  Vanishing theorems}, J. f. Reine
  Angew. Math. 335 (1982), 1-8

\end{thebibliography}

\begin{thebibliography}{EGA}

\bibitem[B]{barton}
C.~M. Barton, \emph{Tensor products of ample vector bundles in
  characteristic $p$}, Amer. J. Math. \textbf{93} (1971), 429--438.
  
\bibitem[G]{gieseker-app}
D. Gieseker, \emph{$p$-ample bundles and their {C}hern classes}, Nagoya
  Math. J. \textbf{43} (1971), 91--116.
  
\bibitem[EGA]{ega-app} 
A.~Grothendieck, \emph{\'{E}l\'ements de g\'eom\'etrie alg\'ebrique}, Inst.
  Hautes \'Etudes Sci. Publ. Math. (1961, 1966), no.~8, 11, 28.
  
\bibitem[H1]{hartshorne-app} 
R. Hartshorne, \emph{Ample vector bundles}, Inst. Hautes \'Etudes Sci. Publ.
  Math. (1966), no.~29, 63--94.

\bibitem[H3]{hartshorne2-app}
\bysame, \emph{Algebraic geometry}, Graduate Texts in Math., no.~52,
  Springer-Verlag, New York, 1977.

\bibitem[Ke]{keeler}
D.~S. Keeler, \emph{Ample filters of invertible sheaves}, {arXiv:math.AG/0108068}, J. Algebra,
 to appear, 2001.
 
\bibitem[Kl]{kleiman}
S.~L. Kleiman, \emph{Toward a numerical theory of ampleness}, Ann. of Math.
  (2) \textbf{84} (1966), 293--344.
 
\end{thebibliography}

\providecommand{\bysame}{\leavevmode\hbox to3em{\hrulefill}\thinspace}
\providecommand{\MR}{\relax\ifhmode\unskip\space\fi MR }
\providecommand{\MRhref}[2]{%
  \href{http://www.ams.org/mathscinet-getitem?mr=#1}{#2}
}
\providecommand{\href}[2]{#2}

\end{appendix}

\end{document}